\documentclass[12pt]{article}

\usepackage{amssymb,amsmath,amsthm}

\newcommand{\ovl}{\overline}

\newcommand{\vp}{\varepsilon}
\newcommand{\cl}[1]{{\mathcal{#1}}}
\newcommand{\bb}[1]{{\mathbb{#1}}}

\newcommand{\A}{{\mathcal A}}
\newcommand{\B}{{\mathcal B}}
\newcommand{\N}{{\mathcal N}}
\newcommand{\LN}{L^2(\cl N,\tr)}
\newcommand{\NB}{\langle \cl N,{\cl B}\rangle}
\newcommand{\tr}{\tau}
\newcommand{\aqual}{\doteq}
\newcommand{\leta}{\ell_{\eta}}
\newcommand{\Leta}{L_{\eta}}

\newcommand{\letaJ}{\ell_{J\eta}}
\newcommand{\NeB}{\N e_{\B}\N}

\numberwithin{equation}{section}

\theoremstyle{plain}
\newtheorem{lem}{Lemma}[section]
\newtheorem{pro}[lem]{Proposition}
\newtheorem{thm}[lem]{Theorem}
\newtheorem{cor}[lem]{Corollary}

\theoremstyle{definition}

\theoremstyle{remark}
\newtheorem{rem}[lem]{Remark}

\begin{document}

\begin{center}
\LARGE  PERTURBATIONS OF SUBALGEBRAS\\
OF TYPE ${\rm {II}}_1$ FACTORS
\end{center}\bigskip

\begin{tabular}{ccc}
Sorin Popa$^*$&Allan M. Sinclair&Roger R.~Smith$^*$ \\
\noalign{\medskip}
Department of Mathematics& Department of Mathematics&Department of 
Mathematics\\
UCLA&University of Edinburgh&Texas A\&M University\\
Los Angeles, CA 90024&Edinburgh, EH9 3JZ&College Station, TX \ 77843\\
U.S.A&SCOTLAND&U.S.A.\\
e-mail:\ {\tt popa@math.ucla.edu}&{\tt A.Sinclair@ed.ac.uk}&{\tt 
rsmith@math.tamu.edu}
\end{tabular}
\vspace{.5in}

\abstract{In this paper we consider two von Neumann subalgebras $\cl B_0$ and 
$\cl
B$ of a type
${\rm{II}}_1$ factor $\cl N$. For a map $\phi$ on $\cl N$, we define
\[\|\phi \|_{\infty,2}=\sup\{\|\phi(x)\|_2\colon \|x\| \leq 1\},\]
and we measure the distance between $\cl B_0$ and $\cl B$ by the quantity
$\|{\bb E}_{\cl B_0}-{\bb E}_{\cl B}\|_{\infty,2}$. Under the hypothesis
that the relative commutant in $\cl N$ of each algebra is equal to its
center,
we prove that close subalgebras have large compressions which are
spatially isomorphic by a partial isometry close to 1 in the $\|\cdot
\|_2$--norm. This hypothesis is satisfied, in particular, by
masas and subfactors of trivial relative commutant.
A general version with a slightly weaker conclusion is also
proved. As a consequence, we show that if $\cl A$ is a masa and
$u\in\cl N$ is a unitary such that $\cl A$ and $u\cl Au^*$ are close, then
$u$ must be close to a unitary which normalizes $\cl A$. These qualitative
statements are given quantitative formulations in the paper.}

\vfill

$\underline{\hspace{1in}}$

\noindent $^*$Partially supported by  grants from the National Science
Foundation.
\newpage

\section{Introduction}\label{sec1}

\indent
In this paper we study pairs of von~Neumann subalgebras $\cl A$ and $\cl B$ of 
a type $\rm{II}_1$ factor $\cl N$ under the assumption that they are close to 
one another in a sense made precise below. Some of our results are very 
general, but the motivating examples are masas, subfactors, or algebras whose 
relative commutants in $\cl N$ equal their centers. In these special cases 
significant extra information is available beyond the general case. Ideally, 
two close subalgebras would be unitarily conjugate by a unitary close to the 
identity, 
but this is not true. In broard terms, we show that two close subalgebras can 
be cut by projections of large trace in such a way that the resulting algebras 
are spatially isomorphic by a partial isometry close to the identity. The 
exact nature of the projections and partial isometry depends on additional 
hypotheses placed on the subalgebras. Our results are an outgrowth of some 
recent work of
the first author who proved a technical rigidity result for two masas
$\cl A$ and $\cl B$ in a type ${\rm {II}}_1$ factor $\cl N$ that has yielded 
several
important results about type ${\rm {II}}_1$ factors, \cite{P3,P4}. The 
techniques of these papers were first developed in \cite{PoCor}.
This paper uses further refinements of these methods to prove 
the 
corresponding stability of certain subalgebras in separable type ${\rm 
{II}}_1$ factors 
(Theorems 
\ref{thm5.3}
and \ref{thm5.4}). Several of the lemmas used below are modifications
of those in \cite{Ch,P3,P4}, and versions of these lemmas go
back to the the foundations of the subject in the papers of 
Murray, von~Neumann, McDuff and Connes. Although the focus of this paper was
initially the
topic of masas, our results have been stated for general von~Neumann algebras
since the proofs are in a similar spirit. The crucial techniques from 
\cite{P3,P4} are
the use of the pull down map $\Phi: L^1(\langle \cl N,{\cl B}\rangle) \to 
L^1(\cl N)$,
and detailed analyses of projections, partial isometries and module properties.
The contractivity of $\Phi$ in the $\|\cdot\|_1$--norm, \cite{P4}, is
replaced here by a discussion of unbounded operators and related norm estimates
in  
Lemma \ref{lem4.1}.

 If $\phi\colon \ \cl N\to
\cl N$ is a bounded linear map, then $\|\phi\|_{\infty,2}$ denotes the quantity
\begin{equation}\label{eq1.1}
\|\phi\|_{\infty,2} = \sup\{\|\phi(x)\|_2\colon \ \|x\|\le 1\},
\end{equation}
and $\|\bb E_{\cl A} - \bb E_{\cl B}\|_{\infty,2}$ measures the distance
between two subalgebras $\cl A$ and $\cl B$, where $\bb E_{\cl A}$ and $\bb 
E_{\cl
B}$ are the associated trace preserving 
conditional expectations. We regard two subalgebras as close to one another if 
$\|\bb E_{\cl A} - \bb E_{\cl B}\|_{\infty,2}$ is small. 
A related notion is that of $\delta$--containment, introduced in 
\cite{MvN} and studied in 
\cite{Ch}. We say that ${\cl A} \subset_{\delta} {\cl B}$ if, for each $a\in 
\cl A$, $\|a\|\leq 1$, there exists $b\in \cl B$ such that $\|a-b\|_2 \leq 
\delta$. This is equivalent to requiring that $\|(I- \bb E_{\cl B})\bb E_{\cl 
A}\|_{\infty,2}\leq \delta$, and so the condition $\|\bb E_{\cl A} - \bb 
E_{\cl 
B}\|_{\infty,2}\leq \delta$ implies $\delta$--containment in both directions, 
so we will often use the norm inequality in the statement of results (see
Remark \ref{rem5.5}).

A significant portion of the paper is devoted to the study of masas.
Two metric based invariants have been introduced to measure the degree of 
singularity of 
a masa $\cl A$ in a type ${\rm {II}}_1$ factor. In \cite{P0}, the delta 
invariant
$\delta(\cl A)$ was introduced, taking values in $[0,1]$. Motivated by this and
certain examples arising from discrete groups, strong singularity and 
$\alpha$--strong
singularity for masas were defined and investigated in \cite{SS1}.
The singular masas are those
which contain their groups of unitary normalizers, \cite{Di}, and within this
class the notion of a strongly singular masa $\cl A\subseteq \cl N$, 
\cite{SS1}, 
is 
defined by the inequality 
\begin{equation}\label{eq1.2}
\|u - \bb E_{\cl A}(u)\|_2 \le \|\bb E_{\cl A} - \bb E_{u\cl A
u^*}\|_{\infty,2}
\end{equation}
for all unitaries $u\in \cl N$. Such an inequality implies that each
normalizing unitary lies in $\cl A$, so (\ref{eq1.2}) can only hold for
singular masas. We may weaken (\ref{eq1.2}) by inserting a constant $\alpha
\in (0,1]$ on the left hand side, and a masa which satisfies the modified
inequality is called $\alpha$-strongly singular. Recently, \cite{P3},
$\delta(\cl A)$ was shown to be 1 for all singular masas. This result,
\cite[Cor. 2]{P3}, may be stated as follows. If $\cl A$ is a singular masa  in 
a 
type
${\rm {II}}_1$ factor $\cl N$ and $v$ is a partial isometry in $\cl N$ with 
$vv^*$ and $v^*v$ orthogonal projections in $\cl A$, then
\begin{equation}\label{eq1.5}
\|vv^*\|_2={\mathrm{sup}}\,\{\|x-{\bb{E}}_{\cl A}(x)\|_2:\ x\in v{\cl A}v^*,\ 
\|x\|
\leq 1\}.
\end{equation}
This result supports the possibility that all singular masas are strongly 
singular.
Although we have not proved that singularity implies strong singularity,
we have been able to establish the inequality
\begin{equation}\label{eq1.8}
\|u - \bb E_{\cl A}(u)\|_2 \le 90\|\bb E_{\cl A} - \bb E_{u\cl A
u^*}\|_{\infty,2}
\end{equation}
for all singular masas $\cl A$ in separably acting type ${\rm 
{II}}_1$ factors $\cl 
N$. 
The constant 90  emerges from a chain of various estimates; our expectation is 
that it should be possible to replace it with a constant equal to or close to 
1.
The method of proof in \cite{P3} (and of the main technical lemma in 
\cite{P4}) 
uses 
the convexity techniques of Christensen, \cite{Ch}, together with the 
pull--down 
identity 
of $\Phi$ from \cite{PP}, some work by Kadison on center--valued traces, 
\cite{K},
and fine estimates on projections. 
Our main proof (Theorem \ref{thm4.2a}) follows that of \cite{P3}, and also 
requires 
approximation of finite projections in 
$L^{\infty}[0,1]\overline{\otimes}B(H)$. 
These 
are combined with a detailed handling of various inequalities involving 
projections and partial isometries.

There are two simple ways in which masas $\cl A$ and $\cl B$ in a type ${\rm 
{II}}_1$ 
factor 
can be close in the $\|\cdot\|_{\infty,2}$--norm on their conditional 
expectations. If $u$
is a unitary close to $\cl A$ in $\|\cdot\|_2$--norm and $u{\cl A}u^*={\cl 
B}$, 
then
$\|\bb E_{\cl A}-\bb E_{\cl B}\|_{\infty,2}$ is small. 
Secondly, if there is a projection $q$ 
of large trace in $\cl A$ and $\cl B$ with $q{\cl A}=q{\cl B}$, then again
$\|\bb E_{\cl A}-\bb E_{\cl B}\|_{\infty,2}$ is small. In Theorem \ref{thm5.4} 
we show 
that a combination of these two methods is the only way in which $\cl A$ can 
be 
close to 
$\cl B$ in separably acting factors.

The structure of the paper is as follows. Section 2 contains preliminary 
lemmas which include statements of some known results that will be used 
subsequently. The operator $h$ in Proposition \ref{pro2.4} was important for 
Christensen's work in \cite{Ch} and plays an essential role here. Theorem 
\ref{thm2.6} investigates its spectrum to aid in later estimates. The third 
section deals with two close algebras, one of which is contained in the other. 
We present a sequence of lemmas, the purpose of which is to cut the algebras 
by a large projection so that equality results. The fourth section collects 
some more background results preparatory to the next section where it is shown 
that two close subalgebras can be cut so that they become isomorphic by a 
suitably chosen partial isometry. In the final section we focus attention on 
applying these results to masas. One consequence is that if a unitary 
conjugate $u\cl A u^*$ of a masa $\cl A$ is close to the original masa then 
$u$ must be close to a normalizing unitary, and this allows us to present the 
results on strongly singular masas mentioned above.

The crucial estimates are contained in Theorems \ref{thm4.2a}, 
\ref{thm3a.5} and Corollary \ref{cor2.5}. We recommend
reading these three results in the order stated, referring back to ancillary
lemmas and propositions as needed. Corollary~\ref{cor2.5} is essentially due
to Christensen in his pioneering paper \cite{Ch}, but without the norm 
inequalities which we have
included. Two 
of our main results, Theorems \ref{thm3a.5} and \ref{thm4.2a}, generalize 
\cite[A.2]{P4} and use methods from \cite[Section 4]{PoCor}.

Our results are formulated for subalgebras of a finite factor $\cl N$. In
only a few places is this requirement necessary, Theorem \ref{thm3a.7} for
example, and when the statement of a result makes sense for  a
von~Neumann algebra $\cl N$ with a unital faithful normal trace, the same 
proof is
valid. 
\newpage

\section{Preliminaries}\label{sec2}

\indent

Let $\cl N$ be a fixed but arbitrary separably acting type ${\rm {II}}_1$ 
factor with
faithful normalized normal trace $\tau$, and let $\cl B$ be a von~Neumann
subalgebra of $\cl N$. The trace induces an inner product
\begin{equation}\label{eq2.1}
\langle x,y\rangle = \tr(y^*x),\qquad x,y\in\cl N,
\end{equation}
on $\cl N$. Then $L^2(\cl N,\tr)$ is the resulting completion with norm
\begin{equation}\label{eq2.2}
\|x\|_2 = (\tr(x^*x))^{1/2},\qquad x\in\cl N,
\end{equation}
and when $x\in\cl N$ is viewed as a vector in this Hilbert space we will
denote it by $\hat x$.
Several traces will be used in the paper and so we will write
$\|\cdot \|_{2,\rm{tr}}$ when there is possible ambiguity. The unique trace 
preserving conditional expectation
$\bb E_{\cl B}$ of $\cl N$ onto $\cl B$ may be regarded as a projection in
$B(L^2(\cl N), \tr)$, where we denote it by $e_{\cl B}$. Thus
\begin{equation}\label{eq2.3}
e_{\cl B}(\hat x) = \widehat{\bb E_{\cl B}(x)},\qquad x\in\cl N.
\end{equation}
Properties of the trace show that there is a conjugate linear isometry
$J\colon \ L^2(\cl N,\tr)\to L^2(\cl N,\tr)$ defined by
\begin{equation}\label{eq2.4}
J(\hat x) = \widehat{x^*},\qquad x\in\cl N,
\end{equation}
and it is standard that $\cl N$, viewed as an algebra of left multiplication
operators on $L^2(\cl N,\tr)$, has commutant $J\cl NJ$. The von~Neumann algebra
generated by $\cl N$ and $e_{\cl B}$ is denoted by $\langle \cl N,{\cl
B}\rangle$, and has commutant $J\cl BJ$. If $\cl B$ is a maximal
abelian self--adjoint subalgebra (masa) of $\cl N$, then $\langle\cl N,{\cl 
B}\rangle$
is a type ${\rm I}_\infty$ von~Neumann algebra, since its commutant is 
abelian. 
Moreover,
its center is $J\cl BJ$, a masa in $\cl N\,'$ and thus isomorphic to
$L^\infty[0,1]$. The general theory of type ${\rm I}$ von~Neumann algebras,
\cite{KR}, shows that there is then a separable infinite dimensional Hilbert 
space
$H$ so that $\langle\cl N,{\cl B}\rangle$ and $L^\infty[0,1]\overline\otimes
B(H)$ are isomorphic. When appropriate, for $\cl B$ a masa, we will regard an 
element $x\in\langle
\cl N,{\cl B}\rangle$ as a uniformly bounded measurable $B(H)$--valued
function $x(t)$ on [0,1]. Under this identification, the center $J\cl BJ$ of
$\langle\cl N,{\cl B}\rangle$ corresponds to those functions taking values
in $\bb CI$. We denote by ${\rm{tr}}$ the unique semi-finite faithful normal 
trace
on $B(H)$ which assigns the value 1 to each rank 1 projection.

The following lemma (see \cite{PP,P2}) summarizes some of the basic properties
of $\langle \cl N,{\cl B}\rangle$ and $e_{\cl B}$.

\begin{lem}\label{lem2.1}
Let $\cl N$ be a separably acting type ${\rm {II}}_1$ factor with a 
von~Neumann 
subalgebra $\cl B$. Then
\begin{align}
\label{eq2.5}
\text{\rm (i)}\quad &e_{\cl B} x e_{\cl B} = e_{\cl B}\bb E_{\cl B}(x) = \bb
E_{\cl B}(x) e_{\cl B},\qquad x\in\cl N;\\
\label{eq2.6}
\text{\rm (ii)}\quad &e_{\cl B}\langle \cl N,{\cl B}\rangle = \overline{e_{\cl
B}\cl
N}^{\,w},\quad \langle \cl N,{\cl B}\rangle e_{\cl B} = \overline{\cl Ne_{\cl
B}}^{\,w};\\
\label{eq2.7}
\text{\rm (iii)}\quad &\text{if } x\in \cl N\cup J{\cl B}J \text{
and } e_{\cl B}x = 0, \text{ then } x=0;\\
\label{eq2.8}
\text{\rm (iv)}\quad &e_{\cl B}\langle \cl N,{\cl B}\rangle e_{\cl B} = \cl
Be_{\cl B} = e_{\cl B}\cl B;\\
\text{\rm (v)}\quad &\text{there is a faithful normal semi-finite trace {\rm
Tr}  on $\langle\cl N,{\cl B}\rangle$ which satisfies}\nonumber
\end{align}
\begin{equation}\label{eq2.9}
{\rm {Tr}}(xe_{\cl B}y) = {\tau}(xy),\qquad x,y\in\cl N,
\end{equation}
and in particular,
\begin{equation}\label{eq2.10}
\text{\rm Tr}(e_{\cl B}) = 1.
\end{equation}
\end{lem}

The following result will be needed subsequently. We denote by
$\|x\|_{2,\text{Tr}}$ the Hilbert space norm induced by Tr on the subspace of
$\langle\cl N,{\cl B}\rangle$ consisting of elements satisfying
$\text{Tr}(x^*x) < \infty$.

\begin{lem}\label{lem2.2}
Let $\cl B$ be a masa in $\cl N$ and let
$\vp>0$. If $f\in \langle \cl N,{\cl B}\rangle$ is a projection of
finite trace and
\begin{equation}\label{eq2.11}
\|f-e_{\cl B}\|_{2,\text{\rm Tr}} \le \vp,
\end{equation}
then there exists a central projection $z\in \langle\cl N,{\cl B}\rangle$
such that $zf$ and $ze_{\cl B}$ are equivalent projections in $\langle\cl N,
e_{\cl B}\rangle$. Moreover, the following inequalities hold:
\begin{equation}\label{eq2.12}
\|zf-ze_{\cl B}\|_{2,\text{\rm Tr}}, \quad \|ze_{\cl B}-e_{\cl
B}\|_{2,\text{\rm Tr}},\quad \|zf-e_{\cl B}\|_{2,\text{\rm Tr}} \le \vp.
\end{equation}
\end{lem}

\begin{proof}
From (\ref{eq2.8}), $e_{\cl B}$ is an abelian projection in $\langle\cl N,
e_{\cl B}\rangle$ so, altering $e_{\cl B}(t)$ on a null set if necessary, each
$e_{\cl B}(t)$ is a projection in $B(H)$ whose rank is at most 1. If $\{t\in
[0,1]\colon \ e_{\cl B}(t) = 0\}$ were not a null set, then there would exist
a non--zero central projection $p \in J{\cl B}J$ corresponding to this set so 
that $e_{\cl
B}p = 0$, contradicting (\ref{eq2.7}). Thus we may assume that each $e_{\cl
B}(t)$ has rank 1.

Since Tr  is a faithful normal semi-finite trace on $\langle\cl N, {\cl
B}\rangle$, there exists a non-negative $\bb R$--valued measurable function
$k(t)$ on [0,1] such that
\begin{equation}\label{eq2.13}
\text{Tr}(y) = \int^1_0 k(t) {\rm{tr}}(y(t))dt, \quad y\in \langle \cl N,{\cl
B}\rangle, \quad \text{Tr}(y^*y) < \infty.
\end{equation}
By (\ref{eq2.10}),
\begin{equation}\label{eq2.14}
\text{Tr}(e_{\cl B}) = \int^1_0 k(t) {\rm{tr}}(e_{\cl B}(t)) dt = 1,
\end{equation}
and thus integration against $k(t)$ defines a probability measure $\mu$ on
[0,1] such that
\begin{equation}\label{eq2.15}
\text{Tr}(y) = \int^1_0 {\rm{tr}}(y(t)) d\mu(t),\quad y\in \langle \cl N, 
e_{\cl
B}\rangle, \quad \text{Tr}(y^*y) < \infty.
\end{equation}
It follows from (\ref{eq2.15}) that
\begin{equation}\label{eq2.16}
\|y\|^2_{2,\text{Tr}} = \int^1_0 {\rm{tr}}(y(t)^* y(t)) d\mu(t),\qquad y\in
\langle\cl N,{\cl B}\rangle.
\end{equation}

Consider a rank 1 projection $p\in B(H)$ and a projection $q\in B(H)$ of rank
$n\ge 2$. Then 
\begin{equation}\label{eq2.17}
{\rm{tr}}((p-q)^2) = {\rm{tr}}(p+q-2pq)
= {\rm{tr}}(p+q-2pqp)
\ge {\rm{tr}}(p+q-2p)
\ge 1,
\end{equation}
and the same inequality is obvious if $q=0$. Let $G = \{t\in [0,1]\colon \
\text{rank}(f(t))\ne 1\}$. Then, from (\ref{eq2.11}),
\begin{equation}\label{eq2.18}
\vp^2 \ge \|f-e_{\cl B}\|^2_{2,\text{Tr}}
\ge \int_G {\rm{tr}}((f(t)-e_{\cl B}(t))^2) d\mu(t)
\ge \mu(G),
\end{equation}
by (\ref{eq2.17}). Let $z = \chi_{G^c} \otimes I$, a central projection in
$\langle\cl N,{\cl B}\rangle$. Then the ranks of $z(t)f(t)$ and
$z(t)e_{\cl B}(t)$ are simultaneously 0 or 1, and so $zf$ and $ze_{\cl B}$ are
equivalent projections in $\langle\cl N,{\cl B}\rangle$. Then
\begin{equation}\label{eq2.19}
\|ze_{\cl B}-e_{\cl B}\|^2_{2,\text{Tr}} = \int_G {\rm{tr}}(e_{\cl B}(t))
d\mu(t)
= \mu(G)
\le \vp^2,
\end{equation}
from (\ref{eq2.18}), while
\begin{align}
\|zf-e_{\cl B}\|^2_{2,\text{Tr}} &= \int_G {\rm{tr}}(e_{\cl B}(t))d\mu(t) +
\int_{G^c} {\rm{tr}}((f(t) - e_{\cl B}(t))^2) d\mu(t)\nonumber\\
\label{eq2.20}
&\le \|f-e_{\cl B}\|^2_{2,\text{Tr}}
\end{align}
since, on $G$,
\begin{equation}\label{eq2.21}
{\rm{tr}}(e_{\cl B}(t)) = 1 \le {\rm{tr}}((f(t) - e_{\cl B}(t))^2),
\end{equation}
by (\ref{eq2.17}). Finally, 
\begin{equation}\label{eq2.22}
\|zf-ze_{\cl B}\|_{2,\text{Tr}} \le \|z\| \|f-e_{\cl B}\|_{2,\text{Tr}} \le
\vp, 
\end{equation}
completely the proof of (\ref{eq2.12}).
\end{proof}

We now recall some properties of the polar decomposition and some trace norm
inequalities. These may be found in \cite{Co,KR}.\newpage

\begin{lem}\label{lem2.3}
Let $\cl M$ be a von Neumann algebra.
\begin{itemize}
\item[\rm (i)] If $w\in\cl M$ then there exists a partial isometry  $v\in\cl
M$, whose initial and final spaces are respectively the closures of the ranges
of $w^*$ and $w$, satisfying
\begin{equation}\label{eq2.23}
w = v(w^*w)^{1/2} = (ww^*)^{1/2}v.
\end{equation}
\item[\rm (ii)] Suppose that $\cl M$ has a faithful normal semifinite trace
\text{\rm Tr}. If $x\in \cl M$, $0\le x\le 1$, \text{\rm Tr}$(x^*x) < \infty$,
and $f$ is the spectral projection of $x$ corresponding to the interval
{\rm [1/2, 1]}, then 
\begin{equation}\label{eq2.24}
\|e-f\|_{2,\text{\rm Tr}} \le 2\|e-x\|_{2,\text{\rm Tr}}
\end{equation}
for any projection $e\in\cl M$ of finite trace.
\item[\rm (iii)]
Suppose that $\cl M$ has a faithful normal semifinite trace and let $p$ and
$q$ be equivalent finite projections in $\cl M$. Then there exists a partial
isometry $v\in \cl M$ and a unitary $u\in\cl M$ satisfying
\begin{gather}\label{eq2.25}
v^*v = p,\quad vv^* = q,\\
\label{eq2.26}
v|p-q|  = |p-q|v,\\
\label{eq2.27}
|v-p|, |v-q| \le 2^{1/2} |p-q|,\\
\label{eq2.28}
upu^* = q,\quad u|p-q| = |p-q|u,\\
\label{eq2.29}
|1-u| \le 2^{1/2} |p-q|.
\end{gather}
\item[\rm (iv)] 
Suppose that $\cl M$ has a faithful normal semifinite trace and let $p$ and
$q$ be finite projections in $\cl M$. Then the partial isometry $v$ in the 
polar decomposition of $pq$ satisfies
\begin{equation}\label{eq2.29a}
\|p-v\|_{2,\text{\rm Tr}},\,\,\|q-v\|_{2,\text{\rm Tr}}\leq \sqrt{2}
\|p-q\|_{2,\text{\rm Tr}}.
\end{equation}
\end{itemize}
\end{lem}
\medskip

The following result is essentially in \cite{Ch}, and is also used in 
\cite{P3,P4}. We reprove it here since the norm estimates that we obtain will
be crucial for subsequent developments. The operator $h$ below will be 
important at several points and we will refer below to the procedure for 
obtaining it as {\it{averaging}} $e_{\cl B}$ {\it{over}} $\cl A$.

\begin{pro}\label{pro2.4}
Let $\cl A$ and $\cl B$ be von~Neumann subalgebras in a separably acting type 
${\rm {II}}_1$ factor $\cl
N$, and let $\ovl K^w_{\cl A}(e_{\cl B})$ be the weak closure of the set
\begin{equation}\label{eq2.30}
K_{\cl A}(e_{\cl B}) = 
{\mathrm{conv}}\,\{ue_{\cl B}u^*\colon \ u \text{ is a unitary in } \cl
A\}
\end{equation}
in $\langle\cl N,{\cl B}\rangle$. Then $\ovl K^w_{\cl A}(e_{\cl B})$
contains a unique element $h$ of minimal $\|\cdot\|_{2,\text{\rm Tr}}$-norm,
and this element satisfies
\begin{align}
\label{eq2.31}
\text{\rm (i)}\quad &h\in \cl A'\cap \langle\cl N,{\cl B}\rangle,\qquad 0
\le h \le 1;\hspace{2.56in}\\
\label{eq2.32}
\text{\rm (ii)}\quad &1-\text{\rm Tr}(e_{\cl B}h) \le \|(I-\bb E_{\cl B})\bb
E_{\cl A}\|^2_{\infty,2};\\
\label{eq2.33}
\text{\rm (iii)}\quad &\text{\rm Tr}(e_{\cl B}h) = \text{\rm Tr}(h^2);\\
\label{eq2.34}
\text{\rm (iv)}\quad &\|h-e_{\cl B}\|_{2,\text{\rm Tr}} \le \|(I-\bb E_{\cl B})
\bb E_{\cl A}\|_{\infty,2}.
\end{align}
\end{pro}

\begin{proof}
Each $x\in K_{\cl A}(e_{\cl B})$ satisfies $0\le x\le 1$, and so the same is
true for elements of $\ovl K^w_{\cl A}(e_{\cl B})$. Moreover, each $x\in
K_{\cl A}(e_{\cl B})$ has unit trace by Lemma \ref{lem2.1} (v). Let $P$ be the
set of finite trace projections in $\langle\cl N,{\cl B}\rangle$. Then,
for $x\in \langle\cl N,{\cl B}\rangle$, $x\ge 0$,
\begin{equation}\label{eq2.35}
\text{Tr}(x) = \sup\{\text{\rm Tr}(xp)\colon \ p\in P\}.
\end{equation}
If $p\in P$ and $(x_\alpha)$ is a net in $K_{\cl A}(e_{\cl B})$ converging
weakly to $x\in \ovl K^w_{\cl A}(e_{\cl B})$, then
\begin{equation}\label{eq2.36}
\lim_\alpha \text{Tr}(x_\alpha p) = \text{Tr}(xp),
\end{equation}
and it follows from (\ref{eq2.35}) that $\text{Tr}(x) \le 1$. Since $x^2\le x$,
it follows that $x\in L^2(\langle\cl N,{\cl B}\rangle,\text{Tr})$. Since
$\text{span}\{P\}$ is  norm dense in 
$L^2(\langle\cl N,{\cl B}\rangle,\text{Tr})$, 
we conclude from (\ref{eq2.36}) that $(x_\alpha)$ converges
weakly to $x$ in the Hilbert space. Thus $\ovl K^w_{\cl A}(e_{\cl B})$ is
weakly compact in both $\langle\cl N,{\cl B}\rangle$ and $L^2(\langle\cl N,
{\cl B}\rangle,\text{Tr})$, and so norm closed in the latter. Thus there
is a unique element $h\in \ovl K^w_{\cl A}(e_{\cl B})$ of minimal
$\|\cdot\|_{2,\text{Tr}}$ -- norm.

For each unitary $u\in\cl A$, the map $x\mapsto uxu^*$ is a
$\|\cdot\|_{2,\text{Tr}}$ -- norm isometry which leaves $\ovl K^w_{\cl
A}(e_{\cl B})$ invariant. Thus 
\begin{equation}\label{eq2.37}
uhu^* = h,\quad u \text{ unitary in } \cl A,
\end{equation}
by minimality of $h$, so $h\in \cl A' \cap \langle \cl N,{\cl B}\rangle$.
This proves (i).

Consider a unitary $u\in\cl A$. Then, by Lemma \ref{lem2.1},
\begin{align}
1 - \text{Tr}(e_{\cl B} ue_{\cl B}u^*) &= 1-\text{Tr}(e_{\cl B}\bb E_{\cl
B}(u)u^*)
= 1-\tr(\bb E_{\cl B}(u)u^*)\nonumber\\
&= 1-\tr(\bb E_{\cl B}(u) \bb E_{\cl B}(u)^*)
= 1-\|\bb E_{\cl B}(u)\|^2_2\nonumber\\
\label{eq2.38}&= \|(I-\bb E_{\cl B}) (u)\|^2_2
\le \|(I-\bb E_{\cl B}) \bb E_{\cl A}\|^2_{\infty,2}.
\end{align}
This inequality persists when $ue_{\cl B}u^*$ is replaced by elements of
$K_{\cl A}(e_{\cl B})$, so it follows from (\ref{eq2.36}) that
\begin{equation}\label{eq2.39}
1-\text{Tr}(e_{\cl B}h) \le \|(I-\bb E_{\cl B}) \bb E_{\cl A}\|^2_{\infty,2},
\end{equation}
proving (ii).

Since $h\in \cl A'\cap \langle\cl N,{\cl B}\rangle$,
\begin{equation}\label{eq2.40}
\text{Tr}(ue_{\cl B}u^*h) = \text{Tr}(e_{\cl B}u^*hu) = \text{Tr}(e_{\cl B}h)
\end{equation}
for all unitaries $u\in\cl A$. Part (iii) follows from this by taking suitable
convex combinations and a weak limit to replace $ue_{\cl B}u^*$ by $h$ on the
left hand side of (\ref{eq2.40}). Finally,  using (\ref{eq2.39}) and
(\ref{eq2.40}),
\begin{align}
\|h-e_{\cl B}\|^2_{2,\text{Tr}} &= \text{Tr}(h^2 - 2he_{\cl B} + e_{\cl
B})\nonumber\\ 
&= \text{Tr}(e_{\cl B}-he_{\cl B})\nonumber\\
&= 1 - \text{Tr}(he_{\cl B})\nonumber\\
\label{eq2.41}
&\le \|(I-\bb E_{\cl B})\bb E_{\cl  A}\|^2_{\infty,2}
\end{align}
proving (iv).
\end{proof}

For the last two results of this  section, $h$ is the element constructed in 
the
previous proposition.

\begin{cor}\label{cor2.5}
Let $\cl A$ and $\cl B$ be von~Neumann subalgebras of $\cl N$ and let
$f$ be the spectral projection of $h$ corresponding to the interval
\newline {\rm [1/2, 1]}. Then $f\in \cl A'\cap\langle\cl N,{\cl B}\rangle$,
and \begin{equation}\label{eq2.42}
\|e_{\cl B}-f\|_{2,\text{\rm Tr}} \le 2\|(I-\bb E_{\cl B})\bb E_{\cl
A}\|_{\infty,2}. \end{equation}
\end{cor}

\begin{proof}
The first assertion is a consequence of elementary spectral theory. The second
follows from Proposition~\ref{pro2.4} (iv) and Lemma~\ref{lem2.3} (ii).
\end{proof}

\begin{thm}\label{thm2.6} 
Let ${\cl Q}_0 \subseteq {\cl Q}_1$ be a containment of finite von~Neumann 
algebras and let $\tau$ be a unital faithful normal trace on ${\cl Q}_1$. 
Suppose that 
\begin{equation}\label{eq2.1a}
{\cl Q}_0'\cap {\cl Q}_1=\cl Z({\cl Q}_0)=\cl Z({\cl Q}_1),
\end{equation}
and let $h \in \langle {\cl Q}_1,{\cl Q}_0\rangle$ be the operator obtained 
from averaging $e_{{\cl Q}_0}$ over ${\cl Q}_1$. Then 
\begin{equation}\label{eq2.2a}
h \in \cl Z({\cl Q}_0)={\cl Q}_1'\cap \langle {\cl Q}_1,{\cl Q}_0\rangle
=\cl Z({\cl Q}_1),
\end{equation}
and the spectrum of $h$ lies in the set 
\begin{equation}\label{eq2.3a}
S=\{(4{\cos}^2(\pi /n))^{-1}: n\geq 3\}\cup [0,1/4].
\end{equation}
In particular, the spectrum of $h$ lies in $\{1\}\cup [0,1/2]$, and the 
spectral projection $q_1$ corresponding to $\{1\}$ is the largest central 
projection for which  ${\cl Q}_0 q_1={\cl Q}_1 q_1$.
\end{thm}

\begin{proof}
Let ${\cl Q}_2$ denote $\langle {\cl Q}_1,{\cl Q}_0\rangle$. Then ${\cl 
Q}_2'=J{\cl Q}_0 J$, so ${\cl Q}_1'\cap {\cl Q}_2=J{\cl Q}_1J\cap (J{\cl 
Q}_0J)'$, which is $J({\cl Q}_0'\cap {\cl Q}_1)J=J(\cl Z({\cl Q}_1))J=\cl 
Z({\cl Q}_1)$. In addition, $\cl Z({\cl Q}_2)=\cl Z({\cl Q}_2')=\cl Z(J{\cl 
Q}_0J)=Z({\cl Q}_1)$, so the algebras $\cl Z({\cl Q}_0)$, $\cl Z({\cl Q}_1)$,
$\cl Z({\cl Q}_2)$, ${\cl Q}_0'\cap {\cl Q}_1$  and ${\cl Q}_1'\cap {\cl Q}_2$
coincide under these hypotheses. Since $h \in {\cl Q}_1'\cap {\cl Q}_2$, by 
Proposition \ref{pro2.4}, we have thus established (\ref{eq2.2a}).

The set $S$ consists of an interval and a decreasing sequence of points. If 
the spectrum of $h$ is not contained in $S$, then we may find a closed 
interval $[a,b]\subseteq S^c\cap(1/4,1)$ so that the corresponding spectral 
projection $z$ of $h$ is non--zero and also lies in  $\cl Z({\cl Q}_1)$. By 
cutting the algebras by $z$, we may assume that $a1\leq h\leq b1$ and 
$[a,b]\cap S=\emptyset$ so $a\geq 1/4$. The trace $\text{Tr}$ on  ${\cl Q}_2$
coming from the basic construction 
satisfies ${\rm {Tr}}(1)\leq a^{-1}{\rm {Tr}}(h)=a^{-1}{\rm {Tr}}(e_{\cl 
Q_0})\leq 4$ and
is thus finite. Let $\text{Ctr}$ denote the 
center--valued trace on ${\cl Q}_2$, whose restrictions to ${\cl Q}_0$ and 
${\cl Q}_1$ are also the  center--valued traces on these subalgebras. Then 
${\text{Ctr}}(e_{{\cl Q}_0})\geq a1$. If ${\cl Q}_2$ has a central summand of 
type ${\rm{I}}_n$ with corresponding central projection $p$ then cutting by 
$p$ gives containment of type ${\rm{I}}_n$ algebras with equal centers and are 
thus equal to each other. This would show that $hp=p$ and $1$ would lie in the 
spectrum of $h$, contrary to assumption. Thus ${\cl Q}_2$ is a type 
$\rm{II}_1$ von~Neumann algebra.
Then $1$ may be expressed as a sum of four equivalent projections 
$\{p_i\}_{i=1}^4$ each having central trace $(4^{-1})1$. Thus each $p_i$ is 
equivalent to a subprojection of $e_{\cl Q_0}$ and so there exist partial 
isometries $v_i \in \cl Q_2$ such that $1=\sum_{i=1}^4 v_ie_{\cl Q_0}v_i^*$. 
Since $e_{\cl Q_0}\cl Q_1=e_{\cl Q_0}\cl Q_2$, we may replace each $v_i$ by an 
operator $w_i \in  \cl Q_1$, yielding
$1=\sum_{i=1}^4 w_ie_{\cl Q_0}w_i^*$. For each $x\in \cl Q_1$, multiply on the 
right by 
$ xe_{\cl Q_0}$ to obtain
\begin{equation}\label{eq2.3b}
xe_{\cl Q_0}=\sum_{i=1}^4 w_i{\bb E}_{{\cl Q_0}}(w_i^*x)e_{\cl Q_0},\ \ x \in 
{\cl Q_1},
\end{equation} 
so $x=\sum_{i=1}^4 w_i{\bb E}_{{\cl Q_0}}(w_i^*x)$, and $\cl Q_1$ is a 
finitely 
generated right $\cl Q_0$--module. In a similar fashion $\cl Q_2=\sum_{i=1}^4 
w_ie_{\cl Q_0}\cl Q_2=\sum_{i=1}^4 w_ie_{\cl Q_0}\cl Q_1$, and so $\cl Q_2$ is 
finitely generated over $\cl Q_1$. This is a standard argument in subfactor 
theory (see \cite{PP}) which we include for the reader's convenience.

Let $\Omega$ be the spectrum of ${\cl Z}(\cl Q_2)$, and fix $\omega \in 
\Omega$. Then
\begin{equation}\label{eq2.5a}
{\cl I}_2=\{x\in \cl Q_2\colon {\text{Ctr}}(x^*x)(\omega)=0\}
\end{equation}
is a maximal norm closed ideal in $\cl Q_2$ and $\cl Q_2/\cl I_2$ is a type 
${\rm{II}}_1$ factor, denoted $\cl M_2$, with trace $\tau_{\omega}=\omega 
\circ {\text{Ctr}}$, \cite{SakY}. Similar constructions yield maximal ideals $
{\cl I}_k\subseteq \cl Q_k$, $k=0,1$, and factors $\cl M_k=\cl Q_k/\cl I_k$. 
Equality of the centers gives $\cl Q_k\cap{\cl I}_2={\cl I}_k$ for $k=0,1$, 
and so $\cl M_0\subseteq \cl M_1\subseteq \cl M_2$ is an inclusion of factors. 
Let $\pi \colon \cl Q_2 \to \cl M_2$ denote the quotient map, and let 
$e=\pi(e_{{\cl Q}_0})$. From above we note that $\cl M_2$ is a finitely 
generated $\cl M_1$--module.

Consider $x \in {\cl I}_1$. By uniqueness of the center--valued trace, the 
composition of ${\bb E}_{{\cl Q_0}}$ with the restriction of Ctr to 
${{\cl Q_0}}$ is Ctr. Thus
\begin{equation}\label{eq2.6a}
{\text{Ctr}}({\bb E}_{{\cl Q_0}}(x^*x)))(\omega)={\text{Ctr}}(x^*x)(\omega)
,\ \ x \in  \cl Q_1.
\end{equation}
Conditional expectations are completely positive unital maps and so
${\bb E}_{{\cl Q_0}}(x^*){\bb E}_{{\cl Q_0}}(x)\leq {\bb E}_{{\cl 
Q_0}}(x^*x)$, showing that $ {\bb E}_{{\cl Q_0}}$ maps $\cl I_1$ to $\cl I_0$.
Thus there is a well defined $\tau_{\omega}$--preserving conditional 
expectation ${\bb E}\colon \cl M_1 \to \cl M_0$ given by 
${\bb E}(x+\cl I_1)={\bb E}_{{\cl Q_0}}(x)+\cl I_0$ for $x \in \cl Q_1$. From 
above, $e$ commutes 
with $\cl M_0$ and $\cl M_2$ is generated by $\cl M_1$ and $e$. Moreover, 
$exe={{\bb E}(x)e}$ for $x \in \cl M_1$ by applying $\pi$ to the equation 
$e_{\cl Q_0}xe_{\cl Q_0}=\bb E_{\cl Q_0}(x)e_{\cl Q_0}$ for $x \in {\cl Q_1}$.
Thus $\cl M_2$ is the extension of $\cl M_1$ by  $\cl M_0$ with Jones 
projection $e$, \cite{Jon}. Now ${\text{Ctr}}(e_{\cl 
Q_0})(\omega)=h(\omega)\in [a,b]$, so $\tau_{\omega}(e)=h(\omega)$ while
$\tau_{\omega}(1)=1$. It follows that $[\cl M_1\colon \cl M_0]^{-1}\in [a,b]$,
contradicting the theorem of Jones, \cite{Jon}, on the possible values of the 
index.

Now let $q_1$ be the spectral projection of $h$ corresponding to $\{1\}$. If 
we cut by $q_1$ then we may assume that $h=1$. But then $e_{\cl Q_0}=1$ and 
$\cl Q_0=\cl Q_1$. We conclude that $\cl Q_0q_1=\cl Q_1q_1$. On the other 
hand, let $z\in \cl Z(\cl Q_1)$ be a projection such that 
$\cl Q_0z=\cl Q_1z$. Then $ e_{\cl Q_0}z=z$, so $hz=z$, showing that $z\leq 
q_1$. Thus $q_1$ is the largest central projection with the stated property. 
\end{proof} 
\newpage

\section{Containment of finite algebras}\label{sec3a}

\indent

In this section we consider an inclusion $\cl M\subseteq\cl N$ of finite von 
Neumann algebras where $\cl N$ has a faithful normal unital trace $\tau$, and 
where $\cl N\subset_\delta\cl M$ for some small positive number $\delta$. Our 
objective is to show that, by cutting the algebras by a suitable projection 
$p$ 
in the center of the relative commutant $\cl M'\cap \cl N$ of large trace 
dependent on $\delta$, we may arrive at $\cl M p = p\cl N p$. This is achieved 
in the following lemmas which are independent of one another. However, we have 
chosen the notation so that they may be applied sequentially to our original 
inclusion $\cl M\subseteq \cl N$. The definition of $\cl N \subset_\delta \cl 
M$ 
depends implicitly on the trace $\tau$ assigned to $\cl N$. Since we will be 
rescaling traces at various points, we will make this explicit by adopting the 
notation $\cl N \subset_{\delta,\tau}\cl M$. If $\tau_1 = \lambda\tau$ for 
some 
$\lambda>0$, then
\begin{equation}\label{eq3a.1}
\|x\|_{2,\tau_1} = \sqrt\lambda \|x\|_{2,\tau},\qquad x\in\cl N.
\end{equation}
Consequently $\cl N\subset_{\delta,\tau} \cl M$ becomes $\cl N 
\subset_{\sqrt\lambda\ \delta,\tau_1} \cl M$ for this change of trace.

It is worth noting that the beginning of the proof of the next lemma shows 
that if $\cl M_0\subseteq  \cl N_0$, $\cl N_0 \subset_{\delta_0,\tau_0} \cl 
M_0$ for a von Neumann algebra $\cl N_0$ with a faithful normal unital trace 
$\tau_0$, then $\cl M'_0\cap \cl N_0 
\subset_{\delta_0,\tau_0} 
\cl Z(\cl M_0)$.

\begin{lem}\label{lem3a.1}
Let $\delta_1\in (0,1)$ and consider an inclusion $\cl M_1\subseteq  \cl N_1$, 
where $\cl N_1$ has a unital faithful normal trace $\tau_1$ relative to which 
$\cl N_1 \subset_{\delta_1,\tau_1} \cl M_1$. Then there exists a projection 
$p_1\in \cl Z(\cl M'_1 \cap \cl N_1)$ such that $\tau_1(p_1)\ge 1-\delta^2_1$ 
and $(\cl M_1p_1)' \cap p_1\cl N_1p_1$ is abelian. 

Setting $\cl M_2  = \cl 
M_1p_1$, $\cl N_2 = p_1\cl N_1p_1$, and $\tau_2 = \tau_1(p_1)^{-1}\tau_1$, we 
have $\cl M'_2\cap \cl N_2$ is abelian and $\cl N_2 \subset_{\delta_2,\tau_2} 
\cl M_2$, where $\delta_2$ is defined by $\delta^2_2 = 
\delta^2_1(1-\delta^2_1)^{-1}$.
\end{lem}

\begin{proof}
Let $\cl C = \cl M'_1\cap \cl N_1$, which contains $\cl Z(\cl M_1)$. 
If $c \in \cl C$, $\|c\|\leq 1$, we may choose $m\in \cl M_1$ to satisfy 
$\|c-m\|_{2,\tau_1}\leq \delta_1$. Conjugation by unitaries from $\cl M_1$ 
leaves $c$ invariant, so
Dixmier's approximation theorem, \cite{Dix}, shows that there is an element  
$z \in \cl Z(\cl M_1)$ such that $\|c-z\|_{2,\tau_1}\leq \delta_1$.
Thus $\cl C 
\subset_{\delta_1} 
\cl Z(\cl M_1)$. Let $\cl A$ be a maximal abelian subalgebra of $\cl C$ which 
contains $\cl Z(\cl M_1)$, and note that $\cl Z(\cl C)\subseteq \cl A$.  
Choose a  projection $p_1\in\cl Z(\cl C)$, maximal with respect to the 
property 
that $\cl Cp_1$ is abelian. We now construct a unitary $u\in \cl C(1-p_1)$ 
such 
that $\bb E_{\cl A(1-p_1)} (u) = 0$.

The algebra $\cl C(1-p_1)$ may be decomposed as a direct sum
\begin{equation}\label{eq3a.2}
\cl C(1-p_1) = \bigoplus_{k\ge 0} \cl C_k,
\end{equation}
where $\cl C_0$ is type $\rm{II}_1$ and each $\cl C_k$ for $k\ge 1$ has the 
form 
$\bb M_{n_k} \otimes \cl A_k$ for an abelian subalgebra $\cl A_k$ of $\cl A$. 
Let $q_k$, $k\ge 0$, be the identity element of $\cl C_k$. Then $\cl A q_k$, 
$k\ge 1$, contains $\cl A_k$ and is maximal abelian in $\cl C_k$, so has the 
form $\cl D_k\otimes \cl A_k$ for some diagonal algebra $\cl D_k \subseteq \bb 
M_{n_k}$, \cite{K}. Note that the choice of $p_1$ implies $n_k\ge 2$. For 
$k\ge  
1$, let $u_k$ be a unitary in $\bb M_{n_k}$ $(\cong \bb M_{n_k} \otimes 1)$ 
which 
cyclically permutes the basis for $\cl D_k$. For $k=0$, choose two equivalent 
orthogonal projections in $\cl A q_0$ which sum to $q_0$, let $v\in \cl C_0$ 
be 
an implementing partial isometry, and let $u_0 = v+v^*$. Then $u   = 
\sum\limits^\infty_{k=0} u_k$ is a unitary in $\cl C(1-p_1)$ for which $\bb 
E_{\cl A(1-p_1)}(u)  = 0$. Thus $w = p_1+u$ is a unitary in $\cl C$. Then
\begin{align}
\delta_1 &\ge \|w-\bb E_{\cl Z(\cl M_1)}(w)\|_{2,\tau_1} \ge \|w-\bb E_{\cl 
A}(w)\|_{2,\tau_1}\nonumber\\
\label{eq3a.3}
&= \|u-\bb E_{\cl A(1-p_1)} (u)\|_{2,\tau_1} = \|u\|_{2,\tau_1} = 
\|1-p_1\|_{2,\tau_1},
\end{align}
and the inequality $\tau_1(p_1) \ge 1 -\delta^2_1$ follows.

Now let $\cl M_2 = \cl M_1p_1$, $\cl N_2 =p_1\cl N_1p_1$ and $\tau_2 
=\tau_1(p_1)^{-1}\tau_1$. Then $\cl M'_2 \cap\cl N_2 = \cl Cp_1$, which is 
abelian, and $\cl N_2 \subset_{\delta_2,\tau_2}\cl M_2$, where $\delta_2 = 
\delta_1(1-\delta^2_1)^{-1/2}$.
\end{proof}

\begin{lem}\label{lem3a.2}
Let $\delta_2\in (0,2^{-1})$ and consider an inclusion $\cl M_2\subseteq \cl 
N_2$, where $\cl N_2$ has a unital faithful normal trace $\tau_2$ relative to 
which $\cl N_2 \subset_{\delta_2,\tau_2} \cl M_2$. Further suppose that $\cl 
M'_2 \cap\cl N_2$ is abelian. Then there exists a projection $p_2\in \cl 
M_2'\cap \cl N_2$ such that $\tau_2(p_2)\ge 1-4\delta^2_2$ and $(\cl M_2p_2)' 
\cap (p_2\cl 
N_2p_2) = \cl Z(\cl M_2p_2)$. In particular, when $\cl N_2$ is abelian we have 
$\cl 
N_2p_2 = \cl M_2p_2$.

Setting $\cl M_3 = \cl M_2p_2$, $\cl N_3 = 
p_2\cl 
M_2p_2$ and $\tau_3=\tau_2(p_2)^{-1}\tau_2$, we have $\cl M'_3 \cap \cl N_3 = 
\cl 
Z(\cl M_3)$ and $\cl N_3 \subset_{\delta_3,\tau_3} \cl M_3$, where $\delta_3$ 
is 
defined by $\delta^2_3 =\delta^2_2(1-4\delta^2_2)^{-1}$.
\end{lem}

\begin{proof}
Let $\cl A = \cl Z(\cl M_2)$ and let $\cl C = \cl M'_2\cap\cl N_2$, which is 
abelian by hypothesis. It is easy to see that $\cl C \subset_{\delta_2,\tau_2} 
\cl A$, by applying Dixmier's approximation theorem, \cite{Dix}. Consider the 
basic 
construction $\cl A\subseteq \cl C \subseteq \langle \cl C,\cl A\rangle$ with 
canonical trace Tr on $\langle \cl C,\cl A\rangle$ given by Tr$(xe_{\cl A} y) 
= \tau_2(xy)$ for $x,y\in\cl C$. Note that $\cl C$ is maximal abelian in 
$B(L^2(\cl C,\tau_2))$ and thus maximal abelian in $\langle \cl C,\cl 
A\rangle$. 
Following the notation of Proposition \ref{pro2.4}, let $h$ be the element of 
minimal $\|\cdot\|_{2,\text{Tr}}$-norm in $\overline K^w_{\cl C}(e_{\cl A})$ 
and 
recall 
from \eqref{eq2.34} that $\|h-e_{\cl A}\|_{2,\text{Tr}} \le \delta_2$. For 
each 
$\lambda\in (2^{-1},1)$, let $f_\lambda$ be the spectral projection of $h$ for 
the interval $[\lambda,1]$. Since $h\in\cl C' \cap \langle \cl C, \cl A\rangle 
= 
\cl C$, we also have that $f_\lambda\in\cl C$ for $2^{-1} < \lambda < 1$. Fix 
an 
arbitrary $\lambda$ in this interval.

We first show that for every projection $q\le f_\lambda$, the inequality
\begin{equation}\label{eq3a.4}
\bb E_{\cl A}(q) \ge \lambda \text{ supp}(\bb E_{\cl A}(q))
\end{equation}
holds. If not, then there exists a projection $q\le f_\lambda$ and $\vp>0$ so 
that the spectral projection $q_0$ of $\bb E_{\cl A}(q)$ for the interval $[0, 
\lambda-\vp]$ is non--zero. Then
\begin{equation}\label{eq3a.5}
0 \ne \bb E_{\cl A}(qq_0) \le \lambda-\vp.
\end{equation}
From this it follows that $e_{\cl A}(qq_0) e_{\cl A} \le (\lambda-\vp) e_{\cl 
A}$, which implies that $qq_0 e_{\cl A} qq_0 \le (\lambda-\vp) 
qq_0$. (To see this, note that, for any pair of projections $e$ and $f$, the 
inequalities $efe\leq \lambda e$, $\|ef\|\leq \sqrt{\lambda}$, $\|fe\|\leq 
\sqrt{\lambda}$, and $fef\leq \lambda f$ are all equivalent). Averaging this 
inequality over unitaries in $\cl Cqq_0$, which have the 
form $uqq_0$ for unitaries $u\in\cl C$, leads to
\begin{equation}\label{eq3a.6}
hqq_0 \le (\lambda-\vp)qq_0.
\end{equation}
The inequality $hf_\lambda \ge \lambda f_\lambda$ implies that
\begin{equation}\label{eq3a.7}
hqq_0 \ge \lambda qq_0,
\end{equation}
and this contradicts \eqref{eq3a.6}, establishing \eqref{eq3a.4}. 

Now consider two orthogonal projections $q_1$ and $q_2$ in $\cl Cf_\lambda$. 
From \eqref{eq3a.4} we obtain
\begin{align}
1 &\ge \bb E_{\cl A}(q_1+q_2) \ge \lambda(\text{supp } \bb E_{\cl A}(q_1) + 
\text{supp } \bb E_{\cl A}(q_2))\nonumber\\
\label{eq3a.8}
&\ge 2\lambda(\text{supp } \bb E_{\cl A}(q_1)\cdot \text{supp } \bb E_{\cl 
A}(q_2)).
\end{align}
Since $\lambda>2^{-1}$, this forces $\bb E_{\cl A}(q_1)$ and $\bb E_{\cl 
A}(q_2)$ to have disjoint support projections. Whenever a conditional 
expectation of one abelian algebra onto another has the property that $\bb 
E(p)\bb E(q)=0$ for all pairs of orthogonal projections $p$ and $q$, then $\bb 
E$ is the identity. This can be easily seen by considering pairs $p$ and 
$1-p$. In our situation, we conclude that
$\cl A f_\lambda = \cl Cf_\lambda$. Let $p_2$ be the spectral projection of 
$h$ 
for the interval $(2^{-1},1]$. By taking the limit $\lambda\to 2^{-1}+$, we 
obtain $\cl Ap_2 = \cl Cp_2$, and the estimate $\tau_2(p_2)\ge 1-4\delta^2_2$ 
follows by taking limits in the inequality
\begin{align}
(1-\lambda)^2 (1-\tau_2(f_\lambda)) &= \tau_2(((1-\lambda)(1-f_\lambda))^2) 
\leq \tau_2((1-h)^2)\nonumber\\
&= \text{Tr}(e_{\cl A}(1-h)^2)= \|e_{\cl A}(1-h)\|^2_{2,\text{Tr}}\nonumber\\
 \label{eq3a.9}&= \|e_{\cl A}(e_{\cl 
A}-h)\|^2_{2,\text{Tr}}
\le \|e_{\cl A}-h\|^2_{2,\text{Tr}} \le \delta^2_2.
\end{align}

Now let $\cl M_3 = \cl  M_2p_2$, $\cl N_3 = p_2\cl N_2p_2$ and $\tau_3 = 
\tau_2(p_2)^{-1}\tau_2$. Then $\cl Z(\cl M_3) = \cl M'_3 \cap \cl N_3$ and 
$\cl 
N_3 \subset_{\delta_3,\tau_3}\cl M_3$, where $\delta_3 = 
\delta_2(1-4\delta^2_2)^{-1/2}$.
\end{proof}

\begin{lem}\label{lem3a.3}
Let $\delta_3 \in (0,4^{-1})$ and consider an inclusion $\cl M_3 \subseteq \cl 
N_3$, where $\cl N_3$ has a unital faithful  normal trace $\tau_3$ relative to 
which $\cl N_3 \subset_{\delta_3,\tau_3} \cl M_3$. Further suppose that $\cl 
Z(\cl M_3) = \cl M'_3 \cap \cl N_3$. Then there exists a projection $p_3\in 
\cl 
Z(\cl M_3)$ such that $\tau_3(p_3)\ge 1-16\delta^2_3$ and
\begin{equation}\label{eq3a.10}
\cl Z(\cl M_3p_3) = (\cl M_3p_3)' \cap (p_3\cl N_3p_3) = \cl Z(p_3\cl N_3p_3).
\end{equation}
Setting $\cl M_4 = \cl M_3p_3$, $\cl N_4 = p_3\cl N_3p_3$ and $\tau_4 = 
\tau_3(p_3)^{-1} \tau_3$, we have
\begin{equation}\label{eq3a.11}
\cl Z(\cl M_4) = \cl M'_4 \cap\cl N_4 = \cl Z(\cl N_4)
\end{equation}
and $\cl N_4 \subset_{\delta_4,\tau_4} \cl M_4$, where $\delta_4$ is defined 
by 
$\delta^2_4 = \delta^2_3(1-16\delta^2_3)^{-1}$.
\end{lem}

\begin{proof}
Since $\cl Z(\cl N_3) \subseteq \cl M'_3\cap \cl N_3$ we have, by hypothesis, 
that $\cl Z(\cl N_3)\subseteq \cl Z(\cl M_3)$. If $x\in\cl Z(\cl M_3)$, 
$\|x\|\le 1$, and $u$ is a unitary in $\cl N_3$ then choose $m\in\cl M_3$ such 
that $\|u-m\|_{2,\tau_3}\le\delta_3$. It follows that
\begin{equation}\label{eq3a.12}
\|ux-xu\|_{2,\tau_3} = \|(u-m)x-x(u-m)\|_{2,\tau_3} \le 2\delta_3,
\end{equation}
and so $\|uxu^*-x\|_{2,\tau_3}\le 2\delta_3$. Suitable convex combinations of 
terms of the form $uxu^*$ converge in norm to an element of $\cl Z(\cl N_3)$, 
showing that $\cl Z(\cl M_3) \subset_{2\delta_3,\tau_3} \cl Z(\cl N_3)$. Now 
apply Lemma \ref{lem3a.2} to the inclusion $\cl Z(\cl N_3) \subseteq \cl Z(\cl 
M_3)$, taking $\delta_2 = 2\delta_3$. We conclude that there is a projection 
$p_3\in \cl Z(\cl M_3)$ such that $\tau_3(p_3)\ge 1-16\delta^2_3$ and $\cl 
Z(\cl 
N_3)p_3 = \cl Z(\cl M_3)p_3$. 

Now let $\cl M_4 = \cl M_3p_3$, $\cl N_4 = \cl N_3p_3$ and $\tau_4 = 
\tau_3(p_3)^{-1}\tau_3$. Then \eqref{eq3a.11} is satisfied and $\cl N_4 
\subset_{\delta_4,\tau_4} \cl M_4$, where 
$\delta_4=\delta_3(1-16\delta^2_3)^{-1/2}$.
\end{proof}

\begin{lem}\label{lem3a.4}
Let $\delta_4\in(0,2^{-1/2})$ and consider an inclusion $\cl M_4 \subseteq \cl 
N_4$, where $\cl N_4$ has a unital faithful normal trace $\tau_4$ relative to 
which $\cl N_4 \subset_{\delta_4,\tau_4}\cl M_4$. Further suppose that
\begin{equation}\label{eq3a.13}
\cl Z(\cl M_4) = \cl M'_4 \cap \cl N_4 = \cl Z(\cl N_4).
\end{equation}
Then there exists a projection $p_4\in \cl Z(\cl M_4)$ such that $\tau_4(p_4) 
\ge 1-2\delta^2_4$ and $\cl M_4p_4 = \cl N_4p_4$.
\end{lem}

\begin{proof}
Consider the basic construction $\cl M_4 \subseteq\cl N_4 \subseteq \langle 
\cl 
N_4, \cl M_4\rangle$ with associated projection $e_{\cl M_4}$, and let $h\in 
\cl 
N'_4 \cap \langle \cl N_4, \cl M_4\rangle$ be the operator obtained from 
$e_{\cl 
M_4}$ by averaging over the unitary group of $\cl N_4$. By hypothesis, the 
conditions of Theorem~\ref{thm2.6} are met, and so $h\in\cl Z(\cl N_4)$ and 
has 
spectrum 
contained in $\{1\}\cup [0,2^{-1}]$. Let $q\in\cl Z(\cl N_4) = \cl Z(\cl M_4)$ 
be the spectral projection of $h$ for the eigenvalue 1, and note that $h(1-q) 
\le (1-q)/2$. Fix an arbitrary $\vp>0$ and suppose that
\begin{equation}\label{eq3a.14}
\text{Tr}(e_{\cl M_4} u(e_{\cl M_4}(1-q))u^*) \ge (2^{-1}+\vp) \text{ 
Tr}(e_{\cl 
M_4}(1-q))
\end{equation}
for all unitaries $u\in\cl M_4$. Taking the average leads to
\begin{equation}\label{eq3a.15}
\text{Tr}(e_{\cl M_4}h(1-q)) \ge (2^{-1}+\vp) \text{ Tr}(e_{\cl M_4}(1-q)),
\end{equation}
and so $\tau_4(h(1-q)) \ge (2^{-1} +\vp) \tau_4(1-q)$. If $q=1$ then we 
already have $\cl N_4 = \cl M_4$; otherwise the last inequality gives a 
contradiction and so \eqref{eq3a.14} fails for every $\vp>0$. The presence of 
$(1-q)$ in \eqref{eq3a.14} ensures that this inequality fails for a unitary 
$u_\vp\in \cl N_4(1-q)$. Thus
\begin{equation}\label{eq3a.16}
\text{Tr}(e_{\cl M_4}u_\vp(e_{\cl M_4}(1-q))u^*_\vp) < (2^{-1}+\vp) \text{ 
Tr}(e_{\cl 
M_4}(1-q)) 
\end{equation}
for each $\vp>0$. Define a unitary in $\cl N_4$ by $v_\vp = q+u_\vp$. By 
hypothesis,
\begin{align}
\delta^2_4 &\ge \|q+u_\vp - \bb E_{\cl M_4}(q+u_\vp)\|^2_{2,\tau_4}
= \|(I-\bb E_{\cl M_4})(u_\vp)\|^2_{2,\tau_4}\nonumber\\
&= \|u_\vp\|^2_{2,\tau_4} - \|\bb E_{\cl M_4}(u_\vp)\|^2_{2,\tau_4}
= \tau_4(1-q) - \tau_4(\bb E_{\cl M_4}(u_\vp)u^*_\vp)\nonumber\\
&= \tau_4(1-q) - \text{Tr}(e_{\cl M_4} u_\vp e_{\cl 
M_4}(1-q)u^*_\vp)\nonumber\\
\label{eq3a.17}
&\ge \tau_4(1-q) - (2^{-1}+\vp) \tau_4(1-q),
\end{align}
where we have used \eqref{eq3a.16} and the fact that
$q\in\cl Z(\cl M_4) = \cl Z(\langle \cl N_4,\cl M_4\rangle)$. Letting $\vp\to 
0$ 
in \eqref{eq3a.17}, we obtain $\tau_4(q)\ge 1-2\delta^2_4$.

Define $p_4 = q\in\cl Z(\cl M_4)$. The basic construction for $\cl M_4 p_4 
\subseteq \cl N_4p_4$ is obtained from the basic construction for $\cl M_4 
\subseteq \cl N_4$ by cutting by the central projection $q$. Since $hq=q$, it 
follows that $\cl N_4p_4 = \cl M_4p_4$, completing the proof.
\end{proof}

We now summarize these lemmas.

\begin{thm}\label{thm3a.5}
Let $\cl N$ be a von Neumann algebra with a
unital faithful normal trace $\tau$, let $\cl M$ be a von Neumann subalgebra, 
and let $\delta$ be a positive number in the interval $(0,(23)^{-1/2})$. If 
$\cl 
N 
\subset_{\delta,\tau}\cl M$, then there exists a projection $p\in\cl Z(\cl 
M'\cap \cl 
N)$ such that $\tau(p)\ge 1-23\delta^2$ and $\cl Mp = p\cl Np$.
\end{thm}

\begin{proof}
We apply the previous four lemmas successively to cut by projections until the 
desired conclusion is reached. Each projection has trace at least a fixed 
proportion of the trace of the previous one, so the estimates in these lemmas 
combine to give
\begin{equation}\label{eq3a.18}
\tau(p)\ge (1-\delta^2_1)(1-4\delta^2_2) (1-16\delta^2_3) (1-2\delta^2_4),
\end{equation}
where the $\delta_i$'s satisfy the relations
\begin{equation}\label{eq3a.19}
\delta^2_1 = \delta^2,\quad \delta^2_2 =\frac{\delta^2_1}{1-\delta^2_1},\quad 
\delta^2_3 = \frac{\delta^2_2}{1-4\delta^2_2},\quad\delta^2_4 = 
\frac{\delta^2_3}{1-16\delta^2_3}.
\end{equation}
Substitution of \eqref{eq3a.19} into \eqref{eq3a.18} gives $\tau(p)\ge 
1-23\delta^2$.
\end{proof}

\begin{rem}\label{rem3a.6}
The assumption that $\delta<(23)^{-1/2}$ in Theorem \ref{thm3a.5} guarantees 
that 
the $\delta_i$'s in the lemmas fall in the correct ranges. This theorem is 
still 
true, but vacuous, for $\delta\ge (23)^{-1/2}$. The constant 23 can be 
improved 
under additional hypotheses by joining the sequence of lemmas at a later 
point. 
If the inclusion $\cl M\subseteq \cl N$, $\cl N\subset_\delta \cl M$ also 
satisfies the hypotheses of Lemmas \ref{lem3a.2}, \ref{lem3a.3} or 
\ref{lem3a.4} 
then 23 can be replaced respectively by 22, 18 or 2.$\hfill\square$
\end{rem}

For the case when the larger algebra is a factor, the estimate in Theorem 
\ref{thm3a.5} can be considerably improved.

\begin{thm}\label{thm3a.7}
Let $\cl N$ be a type ${\rm {II}}_1$ factor with a
unital faithful normal trace $\tau$, let $\cl M$ be a von Neumann subalgebra, 
and let $\delta$ be a positive number in the interval $(0,(2/5)^{1/2})$. If 
$\cl 
N 
\subset_{\delta,\tau}\cl M$, then
 there exists a projection $p\in \cl M'\cap \cl N$ 
with $\tau(p)\ge 1-\delta^2/2$ such that $\cl Mp = p\cl Np$.
\end{thm}

\begin{proof}
Let $q\in\cl M'\cap \cl N$ be a projection with $\tau(q) \ge 1/2$. Then $1-q$ 
is 
equivalent in $\cl N$ to a projection $e\le q$. Let $v\in\cl N$ be a partial 
isometry such that $vv^* = e$, $v^*v = 1-q$. Let $w = v+v^*\in\cl N$ and note 
that $\|w\|=1$ and $\bb E_{\cl M}(w) = 0$. Then
\begin{equation}\label{eq3a.20}
\|w-\bb E_{\cl M}(w)\|^2_{2,\tau} = \|v+v^*\|^2_2 = 2\tau(1-q),
\end{equation}
so we must have $2\tau(1-q)\le \delta^2$, or $\tau(q)\ge 1-\delta^2/2$. On the 
other hand, if $\tau(q) \le 1/2$ then this argument applies to $1-q$, giving 
$\tau(q) \le \delta^2/2$. Thus the range of the trace on projections in $\cl 
M'\cap \cl N$ is contained in $[0,\delta^2/2]\cup [1-\delta^2/2, 1]$.

 By Zorn's lemma and the normality of the trace, there is a 
projection $p\in\cl M'\cap \cl  N$ which is minimal with respect to the 
property 
of having trace at least $1-\delta^2/2$.  We now show that $p$ is a minimal 
projection 
in 
$\cl M'\cap \cl N$. If not, then $p$ can be written $p_1+p_2$ with 
$\tau(p_1)$, 
$\tau(p_2)>0$. By choice of $p$, we see that $\tau(p_1)$, $\tau(p_2) \leq 
\delta^2/2$. It follows that
\begin{equation}\label{eq3a.21}
1-\delta^2/2 \le\tau(p) = \tau(p_1) + \tau(p_2) \leq \delta^2,
\end{equation}
which contradicts $\delta^2<2/5$. Thus $p$ is minimal in $\cl M'\cap \cl 
N$, so $\cl Mp$ has trivial relative commutant in $p\cl Np$. Let $\tau_1$ be 
the 
normalized trace $\tau(p)^{-1}\tau$  on $p\cl Np$. Then $p\cl Np 
\subset_{\delta_1,\tau_1} \cl Mp$, where $\delta_1=\delta(1-\delta^2/2)^{-1/2} 
< 
2^{-1/2}$.

We have now reached the situation of a subfactor inclusion $\cl P\subseteq \cl 
Q$, $\cl Q\subset_\delta\cl P$ for a fixed $\delta<2^{-1/2}$ and $\cl P' \cap 
\cl Q = 
\bb 
C1$. 
Since
\begin{equation}\label{eq3a.27}
\cl Q'\cap \langle\cl Q,\cl P\rangle = J(\cl P' \cap\cl Q)J = \bb C1,
\end{equation}
the operator $h$ obtained from averaging $e_{\cl P}$ over unitaries in $\cl Q$ 
is 
$\lambda 1$ for some $\lambda>0$. By Proposition \ref{pro2.4} (ii) and (iii),
we have $1-\lambda \leq \delta^2<1/2$ and $\lambda = \lambda^2{\rm{Tr}}(1)$, 
yielding ${\rm{Tr}}(1)=1/ \lambda <2$. Thus $[\cl Q \colon \cl P]<2$, so $\cl 
Q = \cl P$ from \cite{Jon}.
Applying this to $\cl Mp \subseteq p\cl Np$, we conclude equality as desired.
\end{proof}

\newpage

\section{Estimates in the $\|\cdot\|_2$-norm}\label{sec3}

\indent

This section establishes some more technical results which will be needed
subsequently. Throughout $\cl N$ is a finite von~Neumann algebra with a unital 
faithful normal trace $\tau$ and $\cl A$ is a general von~Neumann 
subalgebra.

\begin{lem}\label{lem3.1}
Let $w\in \cl N$ have polar decomposition $w=vk$, where $k = (w^*w)^{1/2}$,
and let $p=v^*v$ and $q=vv^*$ be the initial and final projections of $v$. If
$e\in\cl N$ is a projection satisfying $ew=w$, then
\begin{align}
\label{eq3.1}
\text{\rm (i)}\quad &\|p-k\|_2 \le \|e-w\|_2;\hspace{3.3in}\\
\label{eq3.2}
\text{\rm (ii)}\quad &\|e-q\|_2 \le \|e-w\|_2;\hspace{3.3in}\\
\label{eq3.3}
\text{\rm (iii)}\quad &\|e-v\|_2 \le 2\|e-w\|_2.\hspace{3.3in}
\end{align}
\end{lem}

\begin{proof}
The first inequality is equivalent to
\begin{equation}\label{eq3.4}
\tr(p+k^2 - 2pk) \le \tr(e+w^*w -w-w^*),
\end{equation}
since $ew=w$, and (\ref{eq3.4}) is in turn equivalent to
\begin{equation}\label{eq3.5}
\tr(w+w^*) \le \tr(e-p+2k),
\end{equation}
since $pk=k$ from properties of the polar decomposition. The map $x\mapsto
\tr(k^{1/2}xk^{1/2})$ is a positive linear functional whose norm is $\tr(k)$.
 Thus \begin{equation}\label{eq3.6}
|\tr(w)| = |\tr(w^*)| = |\tr(vk)| = |\tr(k^{1/2}vk^{1/2})| \le \tr(k).
\end{equation}
The range of $e$ contains the range of $w$, so $e\ge q$. Thus
\begin{equation}\label{eq3.7}
\tr(e) \ge \tr(q) = \tr(p),
\end{equation}
and so (\ref{eq3.5}) follows from (\ref{eq3.6}), establishing (i).

The second inequality is equivalent to
\begin{align}
\tr(q-2eq) &\le \tr(w^*w - ew-w^*e)\nonumber\\
\label{eq3.8}
&= \tr(k^2-w-w^*).
\end{align}

Since $eq=q$, this is equivalent to
\begin{equation}\label{eq3.9}
\tr(w+w^*) \le \tr(k^2+q) = \tr(p+k^2).
\end{equation}
From (\ref{eq3.6})
\begin{equation}\label{eq3.10}
\tr(w+w^*)\le 2\ \tr(k) = \tr(p+k^2 - (k-p)^2) \le \tr(p+k^2),
\end{equation}
which establishes (\ref{eq3.9}) and proves (ii).

The last inequality is
\begin{align}
\|e-v\|_2 &\le \|e-vk\|_2 + \|v(k-p)\|_2\nonumber\\
\label{eq3.11}&\le \|e-w\|_2 + \|k-p\|_2
\le 2\|e-w\|_2,
\end{align}
by (i).
\end{proof}

The next result gives some detailed properties of the polar decomposition (see
\cite{Ch}). 

\begin{lem}\label{lem3.2}
Let $\cl A$ be a von Neumann subalgebra of $\cl N$ and let
$\phi\colon \ \cl A\to \cl N$ be a normal $*$-homomorphism. Let $w$ have polar
decomposition
\begin{equation}\label{eq3.12}
w = v(w^*w)^{1/2} = (ww^*)^{1/2}v,
\end{equation}
and let $p = v^*v,\ \ q=vv^*$. If
\begin{equation}\label{eq3.13}
\phi(a) w = wa,\qquad a\in \cl A,
\end{equation}
then
\begin{align}
\label{eq3.14}
\text{\rm (i)}\quad &w^*w\in \cl A' \quad \text{and}\quad ww^*\in \phi(\cl
A)';\hspace{2.5in}\\ 
\label{eq3.15}
\text{\rm (ii)} \quad &\phi(a)v = va\quad \text{and}\quad \phi(a)q = vav^*
\text{ for all}\quad a \in \cl A;\hspace{1in}\\\label{eq3.16}
\text{\rm (iii)}\quad &p\in\cl A'\cap \cl N\quad \text{and}\quad q\in\phi(\cl
A)' \cap \cl N.\hspace{2.2in}
\end{align}
\end{lem}

\begin{proof}
If $a\in\cl A$, then
\begin{align}
w^*wa &= w^*\phi(a)w = (\phi(a^*)w)^*w\nonumber\\
\label{eq3.17}
&= (wa^*)^*w = aw^*w,
\end{align}
and so $w^*w \in\cl A'$. The second statement in (i) has a similar proof.

Let $f$ be the projection onto the closure of the range of $(w^*w)^{1/2}$.
Since $w^* = (w^*w)^{1/2}v^*$, the range of $w^*$ is contained in the range of
$f$, and so $f\ge p$ by Lemma~\ref{lem2.3}(i). For all $x\in\cl N$ and
$a\in\cl A$,
\begin{align}
\phi(a) v(w^*w)^{1/2}x &= \phi(a)wx
= wax\nonumber\\
\label{eq3.18}&= v(w^*w)^{1/2}ax
= va(w^*w)^{1/2}x,
\end{align}
since $(w^*w)^{1/2}\in\cl A'$ by (i). Thus
\begin{equation}\label{eq3.19}
\phi(a) vf = vaf,\qquad a\in\cl A,
\end{equation}
which  reduces to
\begin{equation}\label{eq3.20}
\phi(a)v = va,\qquad a\in\cl A,
\end{equation}
since $f\in VN((w^*w)^{1/2}) \subseteq \cl A'$, and
\begin{equation}\label{eq3.21}
v = vp = vpf = vf.
\end{equation}
This proves the first statement in (ii). The second is immediate from
\begin{equation}\label{eq3.22}
\phi(a)q = \phi(a)vv^* = vav^*,\qquad a\in\cl A.
\end{equation}

The proof of the third part is similar to that of the first, and we omit the
details.
\end{proof}

\newpage

\section{Homomorphisms on subalgebras}\label{sec4}

\indent

In this section we consider two close subalgebras $\cl B_0$ and $\cl B$ of a 
type ${\rm {II}}_1$ factor $\cl N$. Our objective is to cut each algebra by a 
projection of large trace in such a way that the resulting algebras are 
spatially isomorphic by a partial isometry which is close to the identity. 
The proof of this involves unbounded operators on $\LN$, so we begin with a
brief discussion of those operators which will appear below. A general
reference for the basic facts about unbounded operators is 
\cite[section 5.6]{KR}.
                                                                                
When $1 \in \N$ is viewed as an element of $\LN$ we will denote it by $\xi$,
and then $x\xi$ is the vector in $\LN$ corresponding to $x\in \N$. We then
have a dense subspace $\N\xi \subseteq \LN$. For each $\eta \in \LN$, we may
define a linear operator $\leta$ with domain $\N\xi$ by
\begin{equation}\label{eq5.1a}
\leta(x\xi)=Jx^*J\eta,\ \ \ x\in \N.
\end{equation}
If $\eta$ happens to be $y\xi$ for some $y \in \N$, then $\leta$ coincides
with $y$, but in general $\leta$ is unbounded. For $x,y\in \N$, we have
\begin{align}\label{eq5.2a}
\langle \leta x\xi,y\xi \rangle &=\langle Jx^*J\eta, y\xi\rangle 
=\overline{\langle x^*J\eta, Jy\xi\rangle } \nonumber\\
&=\langle JyJ\xi, x^*J\eta\rangle  
=\langle x\xi, Jy^*JJ\eta\rangle,
\end{align}
 and so $\leta$ has a densely defined adjoint which agrees with $\letaJ$ on
$\N\xi$. Thus each $\leta$ is closable, and we denote the closure by $\Leta$.
The operators that we consider will all have domains containing $\N\xi$, so
it will be convenient to adopt the notation $S\aqual T$ to mean that $S$ and
$T$
agree on $\N\xi$. This frees us from having to identify precisely the domain
of each particular operator. We note that all unbounded operators arising in
the next result are affiliated with $\N$, and thus any bounded operators
obtained  from the functional calculus will lie in $\N$.

The following lemma will form part of the proof of Theorem \ref{thm4.2a}.
Looking ahead, we will need to draw certain conclusions from (\ref{eq4.15a});
for reasons of technical simplicity we consider the adjoint of this equation
below.
 
\begin{lem}\label{lem4.1}
Let $\B_0$ and $\B$ be von Neumann subalgebras of $\N$, and let $\A$ be a von
Neumann subalgebra of $\B_0$ whose identity is not assumed to be that of $\B_0
$. Suppose that there exist $W\in \NB$ and a $^*$-homomorphism $\theta :\A \to
\B$ such that $aW=W\theta(a)$ for $a\in \A$ and such that $We_{\B}=W$. Then
there exists
a partial isometry $w\in \N$ with the following properties:
 \begin{itemize}
\item[\rm (i)] $w^*aw=\theta(a)w^*w,\ \ a\in \A$;
\item[\rm (ii)] $\|1-w\|_{2,\tau}\leq 2\|e_{\B}-W\|_{2,\text{Tr}}$;
\item[\rm (iii)] $\|1-p'\|_{2,\tau}\leq \|e_{\B}-W\|_{2,\text{Tr}}$,
 where $p'=w^*w\in \theta(\A)'\cap \N$;
\item[\rm (iv)] if $q_1\in \N$ and $q_2\in \B$ are projections such that
$q_1W=W=Wq_2$ then
$q_1w=w=wq_2$.
\end{itemize}
\end{lem}

\begin{proof}
Let $\eta=W\xi\in \LN$. We will first show that $a\Leta \aqual \Leta \theta(a)$
for $a\in \A$. Since ${\text{span}}\{\NeB\}$ is weakly dense in $\NB$, we may
choose, by the Kaplansky density theorem, a sequence $\{y_n\}_{n=1}^{\infty}$
from $ {\text{span}}\{\NeB\}$ converging to $W$ in the strong$^*$ topology.
Since $W=We_{\B}$, we also have that $y_ne_{\B} \to W$ in this topology, and
each $y_ne_{\B}$ has the form $w_ne_{\B}$ for $w_n\in \N$, since 
$e_{\B}\N e_{\B}=\B e_{\B}$ by Lemma \ref{lem2.1}. We also note that $e_{\B}$
commutes
with each $\theta (a)\in \B$, and that $e_{\B}\xi=\xi$. Then, for $a\in \A$,
$x\in \N$,
\begin{align}\label{eq5.4a}
a\Leta x\xi&= aJx^*J\eta=Jx^*JaWe_{\B}\xi \nonumber\\
&=Jx^*JWe_{\B}\theta(a)\xi=\lim_{n\to \infty}Jx^*Jw_ne_{\B}\theta(a)\xi
\nonumber\\
 &= \lim_{n\to \infty}Jx^*Jw_nJ\theta(a^*)J\xi=
\lim_{n\to \infty}Jx^*JJ\theta(a^*)Jw_ne_{\B}\xi \nonumber\\
&=Jx^*JJ\theta(a^*)JW\xi=\Leta \theta(a)x\xi,
\end{align}
establishing that $a\Leta \aqual \Leta \theta(a)$. Let $T=|\Leta |$, and let
$\Leta =wT$ be the polar decomposition  of $\Leta$, where $w$ is a partial
isometry mapping the closure of the range of $T$ to the closure of the range
of $\Leta$. Then
\begin{equation}\label{eq5.5a}
awT\aqual wT\theta(a),\ \ \ a\in \A.
\end{equation}
Let $p'=w^*w$, the projection onto the closure of the range of $T$. Then
$p'T\aqual T$ and (\ref{eq5.5a}) becomes
\begin{equation}\label{eq5.6a}
w^*awT \aqual T\theta(a),\ \ \ a\in \A.
\end{equation}
For each $n\in \mathbb N$, let $e_n\in \N$ be the spectral projection of $T$
for the interval $[0,n]$. Then each $e_n$ commutes with $T$ and so we may
multiply on both sides of (\ref{eq5.6a}) by $e_n$ to obtain
\begin{equation}\label{eq5.7a}
e_nw^*awe_ne_nTe_n=e_nTe_ne_n\theta(a)e_n,\ \ a\in \A,
\end{equation}
where $e_nTe_n\in \N$ is now a bounded operator. When $a\geq 0$, (\ref
{eq5.7a}) implies that $(e_nTe_n)^2$ commutes with $e_n\theta(a)e_n$, and thus
so also does $e_nTe_n$. It follows that
\begin{equation}\label{eq5.8a}
e_nw^*awTe_n=e_n\theta(a)Te_n,\ \ a\in \A,\ n\geq 1.
\end{equation}
For $m\leq n$, we can multiply on the left by $e_m$ and then let $n\to \infty$
to obtain that
\begin{equation}\label{eq5.9a}
e_mw^*aw\zeta=e_m\theta(a)\zeta,\ \ a\in \A, \ \zeta \in {\text{Ran}}\,T.
\end{equation}
Now let $m\to\infty$ to deduce that $w^*aw$ and $\theta(a)$ agree on ${\text
{Ran}}\,T$, and consequently that
\begin{equation}\label{eq5.10a}
w^*aw=w^*aww^*w=\theta(a)w^*w=\theta(a)p',\ \ a \in \A,
\end{equation}
since $p'=1$ on ${\text{Ran}}\,T$. This establishes (i), and by taking $a\geq 0
$ in (\ref{eq5.10a}), it is clear that $p'\in \theta(\A )'\cap \N$. We now turn
to the norm estimates of (ii) and (iii).
                                                                                
On $[0,\infty)$, define continuous functions $f_n$ for $n\geq 1$ by
$f_n(t)=\chi_{[0,n]}(t)+nt^{-1}\chi_{(n,\infty)}(t)$, and let $h_n=f_n(T)\in
\N$ be the associated operators arising from the functional calculus. These
functions were chosen to have the following properties: they form an
increasing sequence with pointwise limit 1, each $f_n$ dominates a positive
multiple of each $\chi_{[0,m]}$, and each $tf_n(t)$ is a bounded function.
Thus $\{h_n\}_{n=1}^{\infty}$ increases strongly to 1 and  each $Th_n$ is a
bounded operator and thus in $\N$. The range of each $h_n$ contains the range
of each $e_m$, and so $\Leta h_n$ and $\Leta$ have identical closures of
ranges for every $n\geq 1$. Thus $w$ is also the partial isometry in the polar
decomposition $\Leta h_n=wTh_n$, $n\geq 1$. The point of introducing the
$h_n$'s is to reduce to the case of bounded operators where we can now apply
Lemma \ref{lem3.1} to obtain
\begin{equation}\label{eq5.11a}
\|1-w\|_{2,\tau}\leq 2\|1-\Leta h_n\|_{2,\tau},\ \
\|1-p'\|_{2,\tau}\leq \|1-\Leta h_n\|_{2,\tau},
\end{equation}
for all $n\geq 1$. For each $b\in \B$,
\begin{equation}\label{eq5.12a}
\Leta b\xi =Jb^*JW\xi=\lim_{n\to \infty}Jb^*Jw_ne_{\B}\xi 
=\lim_{n\to\infty}w_ne_{\B}Jb^*J\xi=Wb\xi,
\end{equation}
using that $e_{\B}$ commutes with both $b$ and $J$. Thus $\Leta e_{\B}=We_{\B}
=W$ and, since $ww^*\Leta \aqual \Leta$, we also have $ww^*W=W$.
Returning to (\ref{eq5.11a}), we obtain
\begin{align}\label{eq5.13a}
\|1-\Leta h_n\|_{2,\tau}&=\|e_{\B}-\Leta h_n e_{\B}\|_{2,{\text{Tr}}}
=\|e_{\B}-wTh_ne_{\B}\|_{2,{\text{Tr}}} \nonumber \\
&=\|e_{\B}-wh_nTe_{\B}\|_{2,{\text{Tr}}}=
\|e_{\B}-wh_nw^*\Leta e_{\B}\|_{2,{\text{Tr}}}\nonumber\\
&=\|e_{\B}-wh_nw^*We_{\B}\|_{2,{\text{Tr}}}.
\end{align}
Since Tr is normal we may let $n\to \infty$ in this last equation, giving
$\lim_{n\to\infty}\|1-\Leta h_n\|_{2,{\tau}}=\|e_{\B}-W\|_{2,{\text{Tr}}}$.
The inequalities of (ii) and (iii) now follow by letting $n\to \infty$ in (\ref
{eq5.11a}). We now establish (iv).
                                                                                
Let $q_1\in \N$ be such that $q_1W=W$.  For each $x\in \N$,
\begin{equation}\label{eq5.14a}
\Leta                                                                         
(x\xi)=Jx^*JW\xi=Jx^*Jq_1W\xi
=q_1Jx^*JW\xi
\end{equation}
and so $q_1$ is the identity on the closure of the range of $\Leta$ which is
also the range of $w$. Thus $q_1w=w$.
                                                                                
Now suppose that $q_2\in \B$ is such that $Wq_2=W$. By replacing $q_2$ by $1-
q_2$, we may prove the equivalent statement that $wq_2=0$ follows from $Wq_2=0
$. Then
\begin{equation}\label{eq5.15a}
0=Wq_2=We_{\B}q_2=Wq_2e_{\B}
=\Leta q_2e_{\B}=wTq_2e_{\B}.
\end{equation}
Multiply by $e_nw^*$ to obtain $e_nTq_2e_{\B}=0$. Since $e_nT\in \N$, it
follows from Lemma \ref{lem2.1} that $e_nTq_2=0$ for all $n\geq 1$. Then
$q_2Te_n=0$ for all $n\geq 1$, and letting $n$ increase, we find that $q_2$
annihilates the closure of the range of $T$ which is also the range of $w^*$.
Thus $q_2w^*=0$ and so $wq_2=0$, completing the proof.
\end{proof}

The following is the main result of this section. We will also state two 
variants which give improved estimates under stronger hypotheses.

\begin{thm}\label{thm4.2a}
Let $\delta>0$, let $\cl B_0$ and $\cl B$ be von Neumann subalgebras of a type 
${\rm {II}}_1$ 
factor $\cl N$ with unital faithful normal trace $\tau$, and suppose that 
$\|\bb 
E_{\cl B} - \bb E_{\cl B_0}\|_{\infty,2}\le \delta$. Then there exist 
projections $q_0\in \cl B_0$, $q\in \cl B$, $q'_0\in \cl B'_0\cap 
\cl 
N$, 
$q'\in \cl B'\cap \cl N$, $p_0' = q_0q'_0$, $p'=qq'$, and a partial isometry 
$v\in 
\cl N$ such that $vp_0'\cl B_0p_0'v^* = p'\cl Bp'$, $vv^* = p'$, $v^*v = 
p_0'$. Moreover, $v$ can be chosen to satisfy 
$\|1-v\|_{2,\tau} \le 69\delta$, $\|1-p'\|_{2,\tau} \le 35\delta$ and 
$\|1-p_0'\|_{2,\tau} \le 
35\delta$.

Under the additional hypothesis that the relative commutants of $\cl B_0$ and 
$\cl B$ are respectively their centers, the projections may be chosen so 
that $p_0'\in\cl B_0$ and $p'\in\cl B$.
\end{thm}

\begin{proof}
We assume that $\delta < (35)^{-1}$, otherwise we may take $v=0$. Let $e_{\cl 
B}$ 
be the Jones projection for the basic construction $\cl B\subseteq \cl N 
\subseteq \langle \cl N,\cl B\rangle$ and let $h\in \cl B'_0\cap\langle \cl 
N,\cl 
B\rangle$ be the operator obtained from averaging $e_{\cl B}$ over 
the unitary group of $\cl B_0$ 
(see Proposition \ref{pro2.4}). If we denote by $e$ the spectral projection of 
$h$ for the interval [1/2,1], then $e\in\cl B'_0 \cap \langle\cl N,\cl 
B\rangle$ 
and $\|e_{\cl B} - e\|_{2,\text{Tr}} \le 2\delta$, by Corollary \ref{cor2.5}. 
Then 
\begin{equation}\label{eq4.6a}
e\cl B_0 = e\cl B_0e \subseteq e\langle \cl N,\cl B\rangle e.
\end{equation}
Consider $x\in e\langle\cl N,\cl B\rangle e$ with $\|x\| \le 1$. Since $e_{\cl 
B}\langle \cl N,\cl B\rangle e_{\cl B} = \cl Be_{\cl B}$, there exists 
$b\in\cl 
B$, $\|b\|\le 1$, such that $e_{\cl B} xe_{\cl B} = be_{\cl B}$. Then
\begin{align}
\|x-\bb E_{\cl B_0}(b)e\|_{2,\text{Tr}}^2 &= \|(x-\bb E_{\cl 
B_0}(b))e\|^2_{2,\text{Tr}}\nonumber\\
\label{eq4.7a}
&= \|(x-\bb E_{\cl B_0}(b)) ee_{\cl B}\|^2_{2,\text{Tr}} + \|(x-\bb E_{\cl 
B_0}(b)) e(1-e_{\cl B})\|^2_{2,\text{Tr}}
\end{align}
and we estimate these terms separately. For the first, we have
\begin{align}
\|(x-\bb E_{\cl B_0}(b))ee_{\cl B}\|_{2,\text{Tr}} &= \|e(x-\bb E_{\cl 
B_0}(b)) 
e_{\cl B}\|_{2,\text{Tr}}\nonumber\\
&\le \|e(x-\bb E_{\cl B}(b)) e_{\cl B}\|_{2,\text{Tr}} + \|e(\bb E_{\cl B}(b) 
- 
\bb E_{\cl B_0}(b)) e_{\cl B}\|_{2,\text{Tr}}\nonumber\\
&\le \|e(x-e_{\cl B}xe_{\cl B})e_{\cl B}\|_{2,\text{Tr}} + \delta\nonumber\\
&= \|e(e-e_{\cl B})(xe_{\cl B})\|_{2,\text{Tr}}  + \delta\nonumber\\
\label{eq4.8a}
&\le \|e-e_{\cl B}\|_{2,\text{Tr}}+ \delta \le 3\delta.
\end{align}
For the second term in \eqref{eq4.7a}, we have
\begin{align}
\|(x-\bb E_{\cl B_0}(b)) e(1-e_{\cl B})\|^2_{2,\text{Tr}} &= \|(x-\bb E_{\cl 
B_0}(b)) e(e-e_{\cl B})\|^2_{2,\text{Tr}}\nonumber\\
&\le \|x-\bb E_{\cl B_0}(b)\|^2 \|e-e_{\cl B}\|^2_{2,\text{Tr}}\nonumber\\
\label{eq4.9a}
&\le 16\delta^2.
\end{align}
Substituting \eqref{eq4.8a} and \eqref{eq4.9a} into \eqref{eq4.7a} gives
\begin{equation}\label{eq4.10a}
\|x-\bb E_{\cl B_0}(b)e\|^2_{2,\text{Tr}} \le 25\delta^2.
\end{equation}
Hence $e\langle\cl N,\cl B\rangle e \subset_{5\delta,\text{Tr}} \cl B_0e$. 
Since 
$\|e-e_{\cl B}\|^2_{2,\text{Tr}}\le 4\delta^2$, it follows that
\begin{equation}\label{eq4.11a}
1+4\delta^2 \ge \text{Tr}(e) \ge 1-4\delta^2.
\end{equation}
If we now  define a unital trace on $e\langle\cl N,\cl B\rangle e$ by $\tau_1 
= 
\text{Tr}(e)^{-1}\text{ Tr}$, then $e\langle\cl N,\cl B\rangle e 
\subset_{\vp,\tau_1} \cl B_0e$ where $\vp = 5\delta(1-4\delta^2)^{-1/2}$. By 
Theorem~\ref{thm3a.5}, there exists a projection $f\in (\cl B_0e)'\cap 
e\langle\cl N,\cl 
B\rangle e$ with $\tau_1(f) \ge 1-23\vp^2$ such that  $\cl B_0f = f\langle \cl 
N,\cl B\rangle f$ (since $ef=f$). 

Let $V\in \langle\cl N,\cl B\rangle$ be the partial isometry in the polar 
decomposition of $e_{\cl B}f$, so that $e_{\cl B}f = (e_{\cl B}fe_{\cl 
B})^{1/2}V$. The inequality $\|V-e_{\cl 
B}\|_{2,\text{Tr}} \le \sqrt 
2\ 
\|e_{\cl B}-f\|_{2,\text{Tr}}$ is obtained from \cite{Co} or Lemma 
\ref{lem2.3} (iv), and we estimate this last quantity. We have
\begin{align}
\|e_{\cl B}- f\|^2_{2,\text{Tr}} &= \text{Tr}(e_{\cl B} + f-2e_{\cl 
B}f)\nonumber\\
&= \text{Tr}(e_{\cl B} + e-2e_{\cl B} e - (e-f) + 2e_{\cl B}(e-f))\nonumber\\
&= \|e_{\cl B}-e\|^2_{2,\text{Tr}} + \text{Tr}(2e_{\cl B}(e-f) - 
(e-f))\nonumber\\
&\le \|e_{\cl B}-e\|^2_{2,\text{Tr}} + \text{Tr}(e-f)
\le 4\delta^2 + 23\vp^2 \text{Tr}(e)\nonumber\\
\label{eq4.13a}
&\le 4\delta^2 + (23)(25) \delta^2(1+4\delta^2)/(1-4\delta^2)
\end{align}
so $\|V-e_{\cl B}\|_{2,\text{Tr}} \le\sqrt 2\ \delta_1$ where $\delta^2_1$ is 
the last quantity above. Then $VV^* \in e_{\cl B}\langle \cl N,\cl B\rangle 
e_{\cl B} = \cl Be_{\cl B}$, and $V^*V\in f\langle \cl N,\cl B\rangle f  = \cl 
B_0 f$, by the choice of $f$. Thus there exist projections $p_0\in\cl B_0$, 
$p\in \cl B$ so that
\begin{equation}\label{eq4.14a}
V^*V = p_0f,\qquad VV^* = pe_{\cl B}.
\end{equation}
If $z\in\cl Z(\cl B_0)$ is the central projection corresponding to the kernel 
of 
the homomorphism $b_0 \mapsto b_0f$ on $\cl B_0$, then by replacing $p_0$ by 
$p_0(1-z)$, we may assume that $p_0b_0p_0 = 0$ whenever $p_0b_0p_0f = 0$.
Since $p_0(1-z)f=p_0f$, (\ref{eq4.14a}) remains valid
and we note that the following relations (and their adjoints) hold:
\begin{equation}\label{eq5.13d}
V=Vp_0=Vf=pV=e_{\cl B}V.
\end{equation}
 Define
$\Theta \colon p_0\cl B_0p_0 \to \langle \cl N, \cl B \rangle$ by
\begin{equation}\label{eq5.13dd}
\Theta(b_0)=Vb_0V^*,\ \ b_0 \in p_0\cl B_0p_0.
\end{equation}
We will show that $\Theta$ is a *--isomorphism onto $p\cl Bpe_{\cl B}$.
Now $pe_{\cl B}V=V$ from (\ref{eq5.13d}), so the range of $\Theta$ is
contained in $pe_{\cl B}\langle \cl N, \cl B \rangle e_{\cl B}p
=p\cl Bpe_{\cl B}$. 
Since
\begin{equation}\label{eq5.13e}
V^*\Theta(b_0)V=p_0fb_0p_0f,\ \ b_0\in p_0\cl B_0p_0,
\end{equation}
from (\ref{eq4.14a}), the choice of $p_0$ shows that $\Theta$ has trivial
kernel. The map is clearly self--adjoint, and we check that it is a
homomorphism. For $p_0b_0p_0,\,p_0b_1p_0 \in p_0\cl B_0p_0$,
\begin{align}
\Theta(p_0b_0p_0)\Theta(p_0b_1p_0)&=Vp_0b_0p_0V^*Vp_0b_1p_0V^*
=Vp_0b_0p_0fp_0b_1p_0V^*\nonumber\\
\label{eq5.13f}&=Vfp_0b_0p_0b_1p_0V^*
=\Theta(p_0b_0p_0b_1p_0),
\end{align}
using $Vf=V$. Finally we show that $\Theta$ maps onto $p\cl Bpe_{\cl B}$.
Given $b\in \cl B$, let
\begin{equation}\label{eq5.13ff}
x=V^*pbpe_{\cl B}V=p_0V^*pbpe_{\cl B}Vp_0 \in p_0f\langle \cl N, \cl B\rangle 
fp_0
=p_0\cl B_0p_0f.
\end{equation}
Then $x$ has the form $p_0b_0p_0f$ for some $b_0 \in \cl B_0$. Thus
\begin{align}
\Theta(p_0b_0p_0)&=Vp_0b_0p_0V^*=Vp_0b_0p_0fV^*\nonumber\\
&=VxV^*=VV^*pbpe_{\cl B}VV^*\nonumber\\
\label{eq5.13g}
&=pe_{\cl B}pbpe_{\cl B}pe_{\cl B}=pbpe_{\cl B},
\end{align}
and this shows surjectivity.
Thus $\Theta \colon p_0\cl B_0p_0 \to p\cl Bpe_{\cl B}$ is a surjective
*--isomorphism, and so can be expressed as
$\Theta(p_0b_0p_0)=\theta(p_0b_0p_0)e_{\cl B}$ where
$\theta \colon p_0\cl B_0p_0 \to p\cl Bp$ is a surjective *--isomorphism.
From the definitions of these maps,
\begin{equation}\label{eq4.15a}
Vb_0 = \theta(b_0)V,\qquad b_0\in p_0\cl B_0p_0.
\end{equation}
If we take the adjoint of this equation then we are in the situation of Lemma
\ref{lem4.1} with $W=V^*$ and $\A =p_0\cl B_0p_0$. We conclude that there is a
partial isometry $v\in \N$ (the $w^*$ of the previous lemma) such that
\begin{equation}\label{eq5.101a}
vb_0v^* = \theta(b_0)vv^*,\qquad b_0\in p_0\cl B_0p_0.
\end{equation}
Then clearly the projection $p'=vv^*$ commutes with $\theta(p_0\B_0 p_0)=p\B
p$, and it lies under $p$ from Lemma \ref{lem4.1} (iv) since $pV=V$. 
 Thus $p'\in (p\B p)'\cap p\N p$.
                                                                                
Now consider the projection $p_0'=v^*v$. Since $Vp_0=V$, it follows from Lemma
\ref{lem4.1} (iv) that $vp_0=v$, and thus $p_0'$ lies under $p_0$. From
(\ref{eq4.15a}), we
have
\begin{equation}\label{eq9.1}
p_0'b_0p_0'=v^*\theta(b_0)v,\ \ \ b_0\in p_0\B_0p_0,
\end{equation}
and the map $b_0\mapsto v^*\theta(b_0)v$ is a *-homomorphism on $p_0\B_0p_0$.
Thus, for all $b_0\in p_0\B_0p_0$,
\begin{equation}\label{eq9.2}
p_0'b_0(1-p_0')b_0^*p_0'=p_0'b_0b_0^*p_0'-(p_0'b_0p_0')(p_0'b_0^*p_0')=0,
\end{equation}
from which we deduce that $p_0'b_0(1-p_0')=0$ and that $p_0'\in (p_0\B_0p_0)'$.
This shows that $p_0'\in (p_0\B_0p_0)'\cap p_0\N p_0$.

It remains to estimate $\|1-v\|_{2,\tau}$. Since
\begin{equation}\label{eq4.21a}
\|V-e_{\cl B}\|_{2,\text{Tr}} \le \sqrt 2\ \|e_{\cl B}-f\|_{2,\text{Tr}} \le 
\sqrt 2\ \delta_1,
\end{equation}
from \eqref{eq4.13a}, we obtain
\begin{equation}\label{eq9.3}
\|1-v\|_{2,\tau}=\|1-w\|_{2,\tau}\leq 2\|e_{\B}-W\|_{2,\text{Tr}}
=2\|e_{\B}-V\|_{2,\text{Tr}}\leq 2\sqrt{2}\delta_1,
\end{equation}
using Lemma \ref{lem4.1} (ii). The estimate of Lemma \ref{lem4.1} (iii) gives
\begin{equation}\label{eq9.4}
\|1-p'\|_{2,\tau}\leq \|e_{\B}-W\|_{2,\text{Tr}}\leq \sqrt{2}\delta_1,
\end{equation}
while a similar estimate holds for $\|1-p_0'\|_{2,\tau}$ because $p'$ and
$p_0'$ are equivalent projections in $\N$.

From the definition of $\delta_1$ and the requirement that $\delta < 
(35)^{-1}$, 
we see that 
\begin{equation}\label{eq4.23a}
8\delta^2_1/\delta^2 \le (69)^2,
\end{equation}
by evaluating the term $(1+4\delta^2) (1-4\delta^2)^{-1}$ at $\delta = 1/35$. 
The 
estimate $\|1-v\|_{2,\tau} \le 69\delta$ follows. Then
\begin{equation}\label{eq4.24a}
\|1-p'\|_{2,\tau}  \le \sqrt 2\ \delta_1 \le 
(69/2)\delta \le 35\delta
\end{equation}
with a similar estimate for $\|1-p'_0\|_{2,\tau}$. The fact that each 
projection is a product of a projection from the algebra and one from the 
relative commutant is clear from the proof. The last statement of the theorem 
is an immediate consequence of the first part, because now the relative 
commutants are contained in the algebras.
\end{proof}

The estimates in Theorem \ref{thm4.2a}, while general, can be substantially 
improved in special cases. The next result addresses the case of two close 
masas.

\begin{thm}\label{thm4.3a}
Let $\delta>0$, let $\cl B_0$ and $\cl B$ be masas in a type ${\rm {II}}_1$ 
factor $\cl 
N$ 
with unital faithful normal trace $\tau$, and suppose that $\|\bb E_{\cl B} - 
\bb 
E_{\cl B_0}\|_{\infty,2} \le \delta$. Then there exists a partial isometry 
$v\in\cl N$ such that $v^*v = p_0\in \cl B_0$, $vv^* = p\in\cl B$, and $v\cl 
B_0 
v^* = \cl Bp$. Moreover $v$ can be chosen to satisfy $\|1-v\|_{2,\tau}\le 
30\delta$, $\|1-p\|_{2,\tau} \le 15\delta$ and $\|1-p_0\|_{2,\tau} \le 
15\delta$.
\end{thm}

\begin{proof}
We assume that $\delta < (15)^{-1}$, otherwise we may take $v=0$. By averaging 
$e_{\cl B}$ over $\cl B_0$, we see that there is a projection $e_0\in \cl B'_0 
\cap 
\langle\cl N,\cl B\rangle$ satisfying $\|e_0-e_{\cl B}\|_{2,\text{Tr}} \le 
2\delta$. By Lemma \ref{lem2.2}, there exists a central projection $z\in 
\langle 
\cl N,\cl B\rangle$ such that $ze_0$ and $ze_{\cl B}$ are equivalent 
projections 
in $\langle\cl N,\cl B\rangle$, and $\|ze_0-e_{\cl B}\|_{2,\text{Tr}} \le 
2\delta$. Let $e = ze_0 \in \cl B'_0 \cap \langle\cl N,\cl B\rangle$, and 
consider the inclusion $\cl B_0 e \subseteq e\langle \cl N,\cl B\rangle e$. 
Let 
$w\in \langle\cl N,\cl B\rangle$ be a partial isometry such that $e=ww^*$, 
$ze_{\cl B} = w^*w$. Then $wze_{\cl B} = w$, and so
\begin{equation}\label{eq4.25a}
\cl B_0 e \subseteq e\langle\cl N,\cl B\rangle e = wze_{\cl B}w^* \langle \cl 
N,\cl B\rangle 
wze_{\cl B}w^* \subseteq wze_{\cl B}\cl B w^*
\end{equation}
and the latter algebra is abelian. The proof now proceeds exactly as in 
Theorem 
\ref{thm4.2a}, starting from (\ref{eq4.6a}) which corresponds to 
\eqref{eq4.25a}. The only difference is that having an 
abelian inclusion allows us to replace 
the 
constant 23 in \eqref{eq4.13a} and subsequent estimates by 4, using Lemma 
\ref{lem3a.2}. 
This leads to the required estimates on $\|1-v\|_{2,\tau}$, $\|1-p\|_{2,\tau}$ 
and $\|1-p_0\|_{2,\tau}$.
\end{proof}

We now consider the case of two close subfactors of $\cl N$.

\begin{thm}\label{thm4.4a}
Let $\delta>0$, let $\cl B_0$ and $\cl B$ be  subfactors of $\cl N$ and 
suppose 
that\break $\|\bb E_{\cl B} - \bb E_{\cl B_0}\|_{\infty,2} \le \delta$. Then 
there 
exist projections $q_0\in \cl B_0$, $q\in \cl B$, $q'_0\in \cl B'_0\cap \cl 
N$, 
$q'\in \cl B'\cap \cl N$, $p_0 = q_0q'_0$, $p=qq'$, and a partial isometry 
$v\in 
\cl N$ such that $vp_0\cl B_0p_0v^* = p\cl Bp$, $vv^* = p$, $v^*v = p_0$, and
\begin{equation}\label{eq4.26a}
\|1-v\|_{2,\tau} \le 13\delta,\qquad \tau(p) = \tau(p_0) \ge 1-67\delta^2.
\end{equation}
If, in addition, the relative commutants of $\cl B_0$ and $\cl B$ are both 
trivial and $\delta<67^{-1/2}$, then $\cl B$ and $\cl B_0$ are unitarily 
conjugate 
in $\cl N$.
\end{thm}

\begin{proof}
We assume that $\delta<67^{-1/2}$, otherwise take $v=0$. The proof is 
identical 
to 
that of Theorem \ref{thm4.2a} except that we now have an inclusion $e\cl B_0e 
\subseteq e\langle\cl N,\cl B\rangle e$ of factors. Our choice of $\delta$ 
allows us a strict upper bound of $(2/5)^{-1/2}$ on the $\vp$ which appears 
immediately after (\ref{eq4.11a}). Thus the estimate of 
Theorem \ref{thm3a.7} applies, which allows us to replace 23 by 1/2 in 
\eqref{eq4.13a}. 
This gives
\begin{equation}\label{eq4.27a}
8\delta^2_1/\delta^2 \le 145 < 169
\end{equation}
and the estimates of \eqref{eq4.26a} follow.

If the relative commutants are trivial then $p\in\cl B$ and $p_0\in\cl B_0$, 
so 
$v$ implements an isomorphism between $p\cl Bp$ and $p_0\cl B_0p_0$ which then 
easily extends to unitary conjugacy between $\cl B$ and $\cl B_0$.
\end{proof}

Let $\cl R$ be the hyperfinite type ${\rm {II}}_1$ factor, choose a projection 
$p\in\cl 
R$ with $\tau(p)  = 1-\delta$, where $\delta$ is small, and let $\theta$ be an 
isomorphism of $p\cl Rp$ onto $(1-p)\cl R(1-p)$. Let $\cl B_0 = 
\{x+\theta(x)\colon \ x\in p\cl Rp\}$ and let $\cl B$ have a similar 
definition 
but using an isomorphism $\phi$ such that $\theta^{-1}\phi$ is a properly 
outer 
automorphism of $p\cl Rp$. Such an example shows that the projections from the 
relative commutants in Theorem \ref{thm4.4a} cannot be avoided.

These results above suggest that it might be possible to obtain similar 
theorems for one sided inclusions. By this we mean that if $\cl B_0
\subset_\delta \cl B$ then there is a partial isometry which moves some  
compression of $\cl B_0$ 
(preferably large) 
into $\cl B$. However the following shows that this 
cannot be so, even if the two algebras are subfactors with trivial relative 
commutant in 
some factor $\cl M$, and even if we renounce 
 the requirement that the size of the 
compression be large and merely require the compression 
to be non--zero. 
In this respect, note that if there exists a 
non--zero partial isometry $v \in \cl M$ such that 
$v^*v \in \cl B_0$, $vv^* \in \cl B$ and 
$vB_0v^* \subseteq vv^*Bvv^*$, then  there would 
be a unitary $u\in \cl M$ such that $u\cl 
B_0u^* \subseteq \cl B$. It is this that we will contradict, 
by exhibiting II$_1$ subfactors $\cl B_0, \cl B \in \cl M$ 
with trivial relative commutant and $\cl B_0
\subset_\delta \cl B$ for $\delta$ arbitrarily small, but with 
no unitary conjugate of $\cl B_0$ sitting inside $\cl B$. 
The construction, based on \cite{PoM}, is given below.  

By \cite{PoM}, 
for each $\lambda < 1/4$ there 
exists an inclusion of factors 
$(\cl N(\lambda) \subseteq \cl M(\lambda))
=(\cl N \subseteq \cl M)$ with Jones index $\lambda^{-1}
> 4$, trivial relative commutant and  
graph $\Gamma_{\cl N,\cl M}=A_\infty$. 
(Note that in fact by \cite{PoS} one can 
take the ambient factor $\cl M$ to be $\cl L(\mathbb F_\infty)$, 
for all $\lambda < 1/4$.) Let $e_0 \in \cl M$ be a 
projection such that $\bb E_{\cl N}(e_0)=\lambda 1$ and let $\cl N_1 \subseteq 
\cl M$ be a subfactor such that $ \cl N_1 \subseteq \cl N \subseteq \cl M$ is 
the basic construction for $ \cl N_1 \subseteq \cl N$ with Jones projection 
$e_0$. Then choose a subfactor $\cl Q \subseteq \cl M$ such that $(1-e_0)\in 
\cl Q$ and $(1-e_0)\cl Q(1-e_0)=\cl N_1(1-e_0)$. An easy 
computation shows that $\cl Q \subset_{\delta(\lambda)} \cl N_1$, 
where $\delta(\lambda) = 6\lambda -4\lambda^2$. Thus, 
since $\cl N_1 \subseteq \cl N$, we get $\cl Q \subset_{\delta(\lambda)} \cl 
N$ 
as well. 

\begin{pro}\label{pro5.5}
With the above notation, we have  
$\cl Q \subset_{\delta(\lambda)} \cl N$, 
with $\delta(\lambda) = 6\lambda - 4\lambda^2$, 
but there does not exist a unitary 
$v\in \cl M$ such that $v\cl Qv^*
\subseteq \cl N$.
\end{pro}

\begin{proof}
Suppose 
there is a unitary $v\in \cl M$ such that $v\cl Qv^*\subseteq \cl 
N_1$, and let $\cl N_0$ be $v^*\cl Nv$. Then $\cl N_0$ is an intermediate 
factor for $\cl Q$. But the irreducible subfactors 
in the Jones tower of a subfactor with Temperley--Lieb--Jones 
standard lattice do not have intermediate subfactors 
(see, for example, \cite{Bis}), giving a contradiction. 

An 
alternative argument goes as follows. 
The basic construction extension algebra 
$\cl Q'\cap \langle \cl M, \cl Q \rangle$ contains 
the projections $e_{\cl Q}$ and $e_{\cl N_0}$, 
which satisfy $e_{\cl Q}\leq e_{\cl N_0}$. Their traces are respectively 
$\lambda^2/(1-\lambda)^2$ and $\lambda$. But the relative commutant 
$\cl Q' \cap \langle \cl M, \cl Q \rangle$ is  
isomorphic to $\bb C^3$
and, of the three minimal projections, the only two 
traces that are less than $1/2$ are 
\begin{equation}\label{eq5.100}
\frac{\lambda^2}{(1-\lambda)^2} {\text{ and }} \frac{\lambda}{1-\lambda}.
\end{equation} 
Since $\tau(e_{\cl N_0})=\lambda$, the only possibility is to have 
\begin{equation}\label{eq5.101}
\tau(e_{\cl N_0})=\frac{\lambda^2}{(1-\lambda)^2}+\frac{\lambda}{1-\lambda}
=\frac{\lambda}{(1-\lambda)^2}.
\end{equation}
This is, of course, impossible.
\end{proof}

\newpage

\section{Unitary congugates of masas}\label{sec5}

\indent

In this section we apply our previous work on perturbations of subalgebras to 
the particular situation of a masa and a nearby unitary conjugate of it. The 
main result of this section is
Theorem~\ref{thm5.3}. This contains two inequalities which we present
separately. Since we will be working with only one unital trace we simplify 
notation by replacing $\|\cdot\|_{2,\tau}$ by $\|\cdot\|_2$, and we denote by 
$d(x,S)$ the distance in $\|\cdot\|_2$-norm from an
element $x\in\cl N$ to a subset $S\subseteq \cl N$. 

Recall from \cite{P1} that the normalizing groupoid $\cl G(\cl A)$ of a masa
$\cl A$ in $\cl N$ is the set of partial isometries $v\in\cl N$ such that
$vv^*$, $v^*v\in\cl A$, and $v\cl Av^* = \cl Avv^*$. Such a partial isometry
$v$ implements a spatial $*$-isomorphism between $\cl Av^*v$ and $\cl Avv^*$.
By choosing a normal $*$-isomorphism between the abelian algebras $\cl
A(1-v^*v)$ and $\cl A(1-vv^*)$ (both isomorphic to $L^\infty[0,1]$), we obtain
a $*$-automorphism of $\cl A$ satisfying the hypotheses of Lemma~\ref{lem2.1}
of \cite{JP}. It follows that $v$ has the form $pw^*$, where $p$ is a
projection in $\cl A$ and $w\in N(\cl A)$ (this result is originally in
\cite{Dy}). The next result will allow us to relate $\|\bb E_{\cl A}  - \bb
E_{u\cl Au^*}\|_{\infty,2}$ to the distance from $u$ to $N(\cl A)$.

\begin{pro}\label{pro3.4}
Let $\cl A$ be a masa in $\cl N$, let $u\in\cl N$ be a unitary and let
$\vp_1,\vp_2>0$. Suppose that there exists a partial isometry $v\in\cl N$ such
that $v^*v\in\cl A$, $vv^*\in u\cl A u^*$, $v\cl Av^* = u\cl Au^*vv^*$, and
\begin{gather}
\label{eq3.32}
\|v-\bb E_{u\cl Au^*}(v)\|_2\le \vp_1,\\
\label{eq3.33}
\|v\|^2_2 \ge 1-\vp^2_2.
\end{gather}
Then there exists $\tilde u\in N(\cl A)$ such that
\begin{equation}\label{eq3.34}
\|u-\tilde u\|_2 \le 2(\vp_1+\vp_2).
\end{equation}
\end{pro}

\begin{proof}
Let $v_1$ be the partial isometry $u^*v\in\cl N$. From the hypotheses we see
that $v^*_1v_1$, $v_1v^*_1\in \cl A$ and
\begin{equation}\label{eq3.35}
v_1\cl Av^*_1 = u^*v\cl Av^*u = u^*u\cl A u^*vv^* u = \cl Av_1v^*_1,
\end{equation}
and so $v_1\in \cl G(\cl A)$. It follows from \cite{JP} that $v_1=pw^*$ for
some projection $p\in\cl A$ and unitary $w^*\in N(\cl A)$. Thus
\begin{equation}\label{eq3.36}
vw = up.
\end{equation}
From (\ref{eq3.32}), there exists $a\in\cl A$ such that $\|a\|\le 1$ and $\bb
E_{u\cl Au^*}(v) = uau^*$. Since $\cl A$ is abelian, it is isomorphic to
$C(\Omega)$ for some compact Hausdorff space $\Omega$. Writing $b=|a|$, $0\le
b\le 1$, there exists a unitary $s\in\cl A$ such that $a=bs$.

Now (\ref{eq3.33}) and (\ref{eq3.36}) imply that
\begin{equation}\label{eq3.37}
\|p\|^2_2 =  \|v\|^2_2 \ge 1-\vp^2_2,
\end{equation}
and so
\begin{equation}\label{eq3.38}
\|1-p\|_2 = (1-\|p\|^2_2)^{1/2} \le \vp_2.
\end{equation}
It now follows from (\ref{eq3.36}) that 
\begin{equation}\label{eq3.39}
\|v-uw^*\|_2 = \|vw-u\|_2 = \|up-u\|_2 = \|1-p\|_2 \le \vp_2.
\end{equation}
From (\ref{eq3.32}) and (\ref{eq3.39}) we obtain the estimate
\begin{align}
\|1-bsu^*w\|_2 &= \|uw^* - ubsu^*\|_2
\le \|uw^*-v\|_2 + \|v-uau^*\|_2\nonumber\\
\label{eq3.40}
&\le \vp_1+\vp_2.
\end{align}
Let $c = \bb E_{\cl A}(su^*w) \in\cl A$, $\|c\|\le 1$, and apply $\bb E_{\cl
A}$ to (\ref{eq3.40}) to obtain
\begin{equation}\label{eq3.41}
\|1-bc\|_2 \le \vp_1+\vp_2.
\end{equation}
For each $\omega\in\Omega$,
\begin{equation}\label{eq3.42}
|\,1-b(\omega)c(\omega)\,|\ge |\,1-|b(\omega)c(\omega)|\,| \ge 1-b(\omega),
\end{equation}
from which it follows that
\begin{equation}\label{eq3.43}
(1-b)^2 \le (1-bc)(1-bc)^*.
\end{equation}
Apply the trace to (\ref{eq3.43}) and use (\ref{eq3.41}) to reach
\begin{equation}\label{eq3.44}
\|1-b\|_2 \le \vp_1+\vp_2.
\end{equation}
Thus
\begin{equation}\label{eq3.45}
\|a-s\|_2 = \|bs-s\|_2 = \|b-1\|_2 \le \vp_1+\vp_2.
\end{equation}
From (\ref{eq3.32}), (\ref{eq3.44}) and the triangle inequality,
\begin{align}
\|v-usu^*\|_2 &= \|v-ubsu^* + u(b-1)su^*\|_2\nonumber\\
\label{eq3.46}&= \|v-\bb E_{u\cl Au^*}(v) + u(b-1)su^*\|_2
\le 2\vp_1 +\vp_2.
\end{align}
This leads to the estimate
\begin{align}
\|u-ws\|_2 
&= \|su^*w-1\|_2
= \|usu^*w-u\|_2\nonumber\\
&= \|usu^*w-up+u(p-1)\|_2
\le \|usu^*w-up\|_2 + \vp_2\nonumber\\
&= \|usu^*w-vw\|_2 + \vp_2
= \|usu^* - v\|_2 + \vp_2\nonumber\\
\label{eq3.47}
&\le 2(\vp_1+\vp_2),
\end{align}
using (\ref{eq3.38}) and (\ref{eq3.46}). Now define $\tilde u = ws$, which is
in $N(\cl A)$ since $s$ is a unitary in $\cl A$. The last inequality gives
(\ref{eq3.34}).
\end{proof}

The constant 90 in the next theorem is not the best possible. An earlier 
version of the paper used methods more specific to masas and obtained the 
lower estimate 31. This may be viewed at the {\it{Mathematics ArXiv}}, 
OA/0111330.

\begin{thm}\label{thm5.1}
Let $\cl A$ be a masa in a separably acting type ${\rm {II}}_1$ factor $\cl 
N$, and
let $u\in\cl N$ be a unitary. Then
\begin{equation}\label{eq5.1}
d(u,N(\cl A)) \le 90 \|(I-\bb E_{u\cl Au^*})\bb E_{\cl A}\|_{\infty,2} \le
90\|\bb E_{\cl A} - \bb E_{u\cl Au^*}\|_{\infty,2}.
\end{equation}
\end{thm}

\begin{proof}
Define $\vp$ to be $\|(I-\bb E_{u\cl Au^*})\bb E_{\cl A}\|_{\infty,2}$. If
$\vp=0$ then $\bb E_{\cl A} = \bb E_{u\cl Au^*}$ and $u\in N(\cl A)$, so there
is nothing to prove. Thus assume $\vp>0$. Let $\cl B = u\cl Au^*$.

By Proposition \ref{pro2.4}, there exists $h\in \cl A'\cap \langle\cl N,{\cl
B}\rangle$ satisfying
\begin{equation}\label{eq5.2}
\|h-e_{\cl B}\|_{2,\text{Tr}}\le \vp.
\end{equation}
Applying Lemma \ref{lem3.2} (ii), the spectral projection $f$ of $h$
corresponding to the interval [1/2,1] lies in $\cl A'\cap \langle\cl N, e_{\cl
B}\rangle$ and satisfies
\begin{equation}\label{eq5.3}
\|f-e_{\cl B}\|_{2,\text{Tr}}\le 2\vp,
\end{equation}
(see Corollary \ref{cor2.5}).
Theorem \ref{thm4.3a} (with $\delta$ replaced by $\vp$) gives the existence of 
a
partial isometry $v\in\cl N$ satisfying
\begin{gather}
\label{eq5.4}
v^*v\in\cl A,\quad vv^*\in\cl B = u\cl Au^*,\quad v\cl Av^* = \cl Bvv^* = u\cl
Au^*vv^*,\\
\label{eq5.5}
\|v-\bb E_{u\cl Au^*}(v)\|_2 \le 30\vp,\\
\label{eq5.6}
\|v\|^2_2 =\tau(vv^*)\ge 1-(15\vp)^2.
\end{gather}
We may now apply Proposition \ref{pro3.4}, with $\vp_1 = 30\vp$ and $\vp_2 = 
15\vp$, to obtain a normalizing unitary
$\tilde u\in N(\cl A)$ satisfying
\begin{equation}\label{eq5.7}
\|u-\tilde u\|_2 \le 2(30 + 15)\vp=90\vp,
\end{equation}
and this is the first inequality.
The second is simply
\begin{align}
\|(I-\bb E_{u\cl Au^*})\bb E_{\cl A}\|_{\infty,2} &= \|(\bb E_{\cl A} - \bb
E_{u\cl Au^*})\bb E_{\cl A}\|_{\infty,2}\nonumber\\
\label{eq5.8}
&\le \|\bb E_{\cl A}-\bb E_{u\cl Au^*}\|_{\infty,2},
\end{align}
completing the proof.
\end{proof}	

\begin{lem}\label{thm5.2}
If $\cl A$ is a von~Neumann subalgebra of a type ${\rm {II}}_1$ factor $\cl N$ 
and 
$u\in\cl N$ is a
unitary, then
\begin{equation}\label{eq5.9}
\|\bb E_{\cl A}-\bb E_{u\cl Au^*}\|_{\infty,2} \le 4d(u, N(\cl A)).
\end{equation}
\end{lem}

\begin{proof}
Let $v\in N(\cl A)$ and define $w$ to be $uv^*$. Then $w\cl Aw^* = u\cl Au^*$,
so it suffices to estimate $\|\bb E_{\cl A}-\bb E_{w\cl Aw^*}\|_{\infty,2}$.
Let $h=1-w$. Then, for $x\in\cl N$, $\|x\|\le 1$,
\begin{align}
\|\bb E_{\cl A}(x) &- \bb E_{w\cl Aw^*}(x)\|_2 = \|\bb E_{\cl A}(x) - w\bb
E_{\cl A}(w^*xw) w^*\|_2\nonumber\\
&= \|w^*\bb E_{\cl A}(x)w-\bb E_{\cl A}(w^*xw)\|_2\nonumber\\
&\le \|w^*\bb E_{\cl A}(x)w - \bb E_{\cl A}(x)\|_2 + \|\bb E_{\cl A}(x) - \bb
E_{\cl A}(w^*xw)\|_2\nonumber\\
&\le \|\bb E_{\cl A}(x)w - w\bb E_{\cl A}(x)\|_2 + \|x-w^*xw\|_2\nonumber\\
&= \|\bb E_{\cl A}(x)h-h\bb E_{\cl A}(x)\|_2 + \|hx-xh\|_2\nonumber\\
\label{eq5.10}&\le 4\|h\|_2
= 4\|1-uv^*\|_2
= 4\|v-u\|_2.
\end{align}
Taking the infimum of the right hand side of (\ref{eq5.10}) over all $v\in
N(\cl A)$ gives (\ref{eq5.9}).
\end{proof}

The next theorem summarizes the previous two results.

\begin{thm}\label{thm5.3}
Let $\cl A$ be a masa in a separably acting type ${\rm {II}}_1$ factor $\cl N$ 
and let
$u$ be a unitary in $\cl N$. Then
\begin{equation}\label{eq5.11}
d(u,N(\cl A))/90 \le \|(I-\bb E_{u\cl Au^*})\bb E_{\cl A}\|_{\infty,2} \le
\|\bb E_{\cl A}-\bb E_{u\cl Au^*}\|_{\infty,2} \le 4d(u,N(\cl A)).
\end{equation}
If $\cl A$ is singular, then $\cl A$ is (1/90)-strongly singular.
\end{thm}

\begin{proof}
The inequalities of (\ref{eq5.11}) are proved in Theorem \ref{thm5.1} and
Lemma \ref{thm5.2}. When $\cl A$ is singular, its normalizer is contained in 
$\cl
A$, so
\begin{equation}\label{eq5.12}
\|u-\bb E_{\cl A}(u)\|_2 \le d(u, N(\cl A))
\end{equation}
holds. Then
\begin{equation}\label{eq5.13}
\|u-\bb E_{\cl A}(u)\|_2 \le 90\|\bb E_{\cl A} - \bb E_{u\cl
Au^*}\|_{\infty,2},
\end{equation}
proving $\alpha$-strong singularity with $\alpha = 1/90$.
\end{proof}

The right hand inequality of (\ref{eq5.11}) is similar to 
\begin{equation}\label{eq5.14}
\|\bb E_{\cl A} - \bb E_{u\cl
Au^*}\|_{\infty,2}\leq 4\|u-\bb E_{\cl A}(u)\|_2,
\end{equation}
which we obtained in \cite[Prop. 2.1]{SS1}, so $u$ being close to $\cl A$ 
implies
that $\cl A$ and $u{\cl A}u^*$ are also close. We remarked in the introduction 
that there are only two ways in which $\|\bb E_{\cl A} - \bb E_{u\cl
Au^*}\|_{\infty,2}$ can be small, and we now make precise this assertion and 
justify it.

\begin{thm}\label{thm5.4}
Let $\cl A$ and $\cl B$ be masas in a separably acting type ${\rm {II}}_1$ 
factor $\cl 
N$,
and let $\delta_1,\,\delta_2,\,\vp > 0$.
\begin{itemize}
\item[\rm (i)] If there are projections $p \in \cl A$, $q \in \cl B$ and a 
unitary
$u \in \cl N$ satisfying
\begin{equation}\label{eq5.15}
u^*qu=p,\ \ \ u^*q{\cl B}u=p{\cl A},
\end{equation}
\begin{equation}\label{eq5.16}
\|u-{\bb E}_{\cl B}(u)\|_2 \leq \delta_1
\end{equation}
and
\begin{equation}\label{eq5.17}
\text{tr}(p)=\text{tr}(q)\geq 1-{\delta_2}^2,
\end{equation}
then
\begin{equation}\label{eq5.18}
\|\bb E_{\cl A} - \bb E_{\cl
B}\|_{\infty,2}\leq 4\delta_1 + 2\delta_2.
\end{equation}
\item[\rm (ii)] If $\|\bb E_{\cl A} - \bb E_{\cl
B}\|_{\infty,2}\leq \vp$, then there are projections $p \in \cl A$ and $q \in 
\cl B$,
and a unitary $u \in \cl N$ satisfying
\begin{equation}\label{eq5.19}
u^*qu=p,\ \ u^*q{\cl B}u=p\cl A,
\end{equation}
\begin{equation}\label{eq5.20}
\|u-\bb E_{\cl B}(u)\|_2 \leq 45\vp
\end{equation}
and
\begin{equation}\label{eq5.21}
\text{tr}(p)=\text{tr}(q)\geq 1-(15\vp)^2.
\end{equation}

\end{itemize}
\end{thm}

\begin{proof} 
\noindent (i) Let $\cl C=u^*{\cl B}u$. Then 
\begin{equation}\label{eq5.22}
\|\bb E_{\cl B} - \bb E_{\cl
C}\|_{\infty,2}\leq 4\delta_1,
\end{equation}
 from (\ref{eq5.14}). If $x \in \cl N$ with 
$\|x\|\leq 1$, then 
\begin{align}
\|\bb E_{\cl C}(x)-\bb E_{\cl A}(x)\|_2&\leq\|(1-p)(\bb E_{\cl C}-\bb E_{\cl 
A})(x)\|_2
+\|p\bb E_{\cl C}(px)-\bb E_{\cl A}(px)\|_2\nonumber\\
\label{eq5.23}&\leq 2\delta_2,
\end{align}
since
\begin{equation}\label{eq5.24}
p{\cl C}=pu^*{\cl B}u=u^*q{\cl B}u=p{\cl A}
\end{equation}
and $p\bb E_{\cl A}(p(\cdot))$ is the projection onto $p{\cl A}$. Then
(\ref{eq5.18}) follows immediately from (\ref{eq5.22}) and (\ref{eq5.23}).

\noindent (ii) As in the proofs of Theorem \ref{thm5.1} and its preceding 
results
Proposition \ref{pro2.4}, Lemma \ref{lem3.2} and Theorem \ref{thm4.2a},
there is a partial isometry $v \in \cl N$ satisfying
\begin{equation}\label{eq5.25}
p=v^*v \in {\cl A},\ \ q=vv^* \in {\cl B},\ \ v^*q{\cl B}v=p{\cl A},
\end{equation}
\begin{equation}\label{eq5.26}
\|v-\bb E_{\cl B}(v)\|_2 \leq 30\vp
\end{equation}
and
\begin{equation}\label{eq5.27}
\text{tr}(p)=\text{tr}(q)\geq 1 -(15\vp)^2.
\end{equation}
Let $w$ be a partial isometry which implements the equivalence
\begin{equation}\label{eq5.28}
w^*w=1-p,\ \ \ ww^*=1-q,
\end{equation}
and let $u=v+w$. Then $u$ is a unitary in $\cl N$, since the initial
and final projections of $v$ and $w$ are orthogonal, and
\begin{equation}\label{eq5.29}
u^*q{\cl B}u=v^*q{\cl B}v=p {\cl A}.
\end{equation}
Observe that 
\begin{equation}\label{eq5.30}
\|w-\bb E_{\cl B}(w)\|_2 \leq \|w\|_2 =(\text{tr}(1-p))^{1/2}\leq 15\vp,
\end{equation}
so that the inequality
\begin{equation}\label{eq5.31}
\|u-\bb E_{\cl B}(u)\|_2\leq 45\vp
\end{equation}
follows from (\ref{eq5.26}) and (\ref{eq5.30}).
\end{proof}

\begin{rem}\label{rem5.5}
Recall that ${\cl A}\subset_{\delta} {\cl 
B}$ 
is equivalent to $\|(I-\bb E_{\cl B})\bb E_{\cl A}\|_{\infty,2}
\leq \delta$. In \cite{Ch}, Christensen defined the distance between $\cl A$ 
and $\cl B$ to be
\begin{equation}\label{eq5.32}
\|{\cl A}-{\cl B}\|_2=\text{max}\,\{\|(I-\bb E_{\cl B})\bb E_{\cl 
A}\|_{\infty,2},\ 
\|(I-\bb E_{\cl A})\bb E_{\cl B}\|_{\infty,2}\}.
\end{equation}
This quantity is clearly bounded by $\|\bb E_{\cl A}-\bb E_{\cl 
B}\|_{\infty,2}$,
and the reverse inequality 
\begin{equation}\label{eq5.33}
\|\bb E_{\cl A}-\bb E_{\cl B}\|_{\infty,2}\leq 3\|{\cl A}-{\cl B}\|_2
\end{equation}
follows from \cite[Lemma 5.2]{SS1} and the algebraic identity
\begin{equation}\label{eq5.34}
P-Q=P(I-Q)-(I-P)Q,
\end{equation}
valid for all operators $P$ and $Q$. 
In general, for any $x\in \cl N$, $\|x\|\le 1$,
\begin{align}
\|(\bb E_{\cl A}-\bb E_{\cl B})(x)\|^2_2 &= \langle \bb E_{\cl A}(x), (\bb
E_{\cl A}-\bb E_{\cl B})(x)\rangle - \langle \bb E_{\cl B}(x), (\bb E_{\cl A}
- \bb E_{\cl B})(x)\rangle\nonumber\\
&= \langle (I-\bb E_{\cl B})\bb E_{\cl A}(x),x\rangle + \langle(I-\bb E_{\cl
A}) \bb E_{\cl B}(x), x\rangle\nonumber\\
\label{eq5a.23}
&\le \|(I-\bb E_{\cl B})\bb E_{\cl A}\|_{\infty,2} + \|(I-\bb E_{\cl A}) \bb
E_{\cl B}\|_{\infty,2},
\end{align}
which gives the inequality
\begin{equation}\label{eq5a.33}
\|\bb E_{\cl A}-\bb E_{\cl B}\|_{\infty,2}\leq (2\|{\cl A}-{\cl B}\|_2)^{1/2}.
\end{equation}
Thus the two notions of distance give equivalent metrics on the space of all 
subalgebras of $\cl N$. $\hfill\square$ 
\end{rem}

We close with a topological result on the space of masas, in the spirit of 
\cite{Ch,SS1},
which also follows from results in \cite{P3}. We include a short proof for 
completeness.

\begin{cor}\label{cor5.6}
The set of singular masas in a separably acting type ${\rm {II}}_1$ factor is 
closed
in the $\|\cdot\|_{\infty,2}$-metric.
\end{cor}

\begin{proof}
By Theorem \ref{thm5.3}, it suffices to show that those masas, which satisfy
(\ref{eq5.13}) (with any fixed $\alpha>0$ replacing 90) for all unitaries
$u\in \cl N$, form a closed subset. Consider a Cauchy sequence $\{\cl
A_n\}^\infty_{n=1}$ of masas satisfying (\ref{eq5.12}), and fix a unitary
$u\in\cl N$. By \cite{Ch}, the set of masas is closed, so there is a masa $\cl
A$ such that $\lim\limits_{n\to\infty} \|\bb E_{\cl A_n} - \bb E_{\cl
A}\|_{\infty,2} = 0$. Then
\begin{align}
\|u-\bb E_{\cl A}(u)\|_2 &\le \|u-\bb E_{\cl A_n}(u)\|_2 + \|\bb E_{\cl
A_n}(u) - \bb E_{\cl A}(u)\|_2\nonumber\\
&\le \alpha\|\bb E_{u\cl A_nu^*} - \bb E_{\cl A_n}\|_{\infty,2} + \|\bb
E_{\cl A_n} - \bb E_{\cl A}\|_{\infty,2}\nonumber\\
\label{eq5.35}
&\le \alpha\|\bb E_{u\cl A_nu^*} - \bb E_{u\cl Au^*}\|_{\infty,2} + \alpha
\|\bb E_{u\cl Au^*} - \bb E_{\cl A}\|_{\infty,2}
+ \|\bb E_{\cl A_n} - \bb E_{\cl A}\|_{\infty,2},
\end{align}
and the result follows by letting $n\to \infty$.
\end{proof}

\end{document}